\documentclass[11pt]{article}
\usepackage{mathrsfs}

\usepackage{amssymb}
\usepackage{mathrsfs}
\usepackage{amsfonts}
\usepackage{amsmath}
\usepackage{cite}
\usepackage{algorithm}
\usepackage{algorithmic}

\newcommand{\ba}{\noindent $\begin{array}}
\newcommand{\ea}{\end{array}$}
\newcommand{\be}{\begin{equation}}
\newcommand{\ee}{\end{equation}}
\newcommand{\bd}{\begin{displaymath}}
\newcommand{\ed}{\end{displaymath}}
\newcommand{\beq}{\begin{eqnarray*}}
\newcommand{\eeq}{\end{eqnarray*}}
\newcommand{\beqn}{\begin{eqnarray}}
\newcommand{\eeqn}{\end{eqnarray}}

\usepackage{color}

\newtheorem{theorem}{Theorem}[section]
\newtheorem{proposition}{Proposition}[section]
\newtheorem{lemma}{Lemma}[section]
\newtheorem{corollary}{Corollary}[section]
\newtheorem{definition}{Definition}[section]

\newtheorem{remark}{Remark}[section]

\newfont{\Bb}{msbm10 scaled\magstep1}

\def\sqr#1#2{{\vcenter{\hrule height .#2pt
      \hbox{\vrule width .#2pt height#1pt \kern#1pt\vrule width.#2pt}
                       \hrule height.#2pt}}}

\def\argmin{\,{\rm argmin}\,}

\def\K{\mathcal{K}}

\title{Variational Analysis of the Orthogonally Invariant  Norm Cone of Symmetric Matrices
}
\author{Yule Zhang\footnote{School of Science, Dalian Maritime University, Dalian 116026, China. (ylzhang@dlmu.edu.cn) This author was supported by the Natural Science Foundation of China under No. 12201097.},\,\, Jihong Zhang\footnote{School of Science, Shenyang Polytechnic University, Shenyang 110159, China. (zjh7815040x@163.com)}
 \,\,and \,\,Liwei Zhang
\footnote{School of Mathematical Sciences, Dalian University of Technology, Dalian 116024, China. (lwzhang@dlut.edu.cn) This author was supported by National Key R\&D Program of China under project No. 2022YFA1004000, the Natural Science Foundation of China under No. 11971089 and partially supported by Dalian High-level Talent Innovation Project No. 2020RD09.}
}

\date{\today}

\begin{document}
\maketitle

 \vspace{2mm}
\begin{center}
\parbox{13.5cm}{\small
\textbf{Abstract.} A large number matrix optimization problems  are described by
  orthogonally invariant  norms. This paper is devoted to the study of variational analysis of the orthogonally invariant  norm cone of symmetric matrices. For a general orthogonally invariant norm cone of symmetric matrices, formulas for the tangent cone, normal cone  and second-order tangent set are established. The differentiability properties of the projection operator onto the orthogonally invariant norm cone are developed, including formulas for the directional derivative and the B-subdifferential. Importantly, the directional derivative is characterized by the second-order derivative of the corresponding symmetric function, which is convenient for computation. All these results are specified to the Schatten $p$-norm cone, especially to  the second-order cone of symmetric matrices. \\

 }
\end{center}

\section{\textbf{Introduction}}
\setcounter{equation}{0}
We consider the  space $\mathbb S^m$ of real $m\times m$ symmetric matrices endowed with standard inner product
$\langle \cdot, \cdot \rangle$ by $\langle A, B\rangle={\rm Tr}, (A^TB)$ for $A,B \in \mathbb S^m$. The induced norm is the Frobenius  norm, denoted by $\|\cdot\|_F$:
$\|A\|_F=\sqrt{\langle A, A\rangle}$ for $ A \in \mathbb S^m$.
We say that a norm $\|\cdot\|$ in $\mathbb S^m$ is an orthogonal invariant norm  if
$$
\|P^TAP\|=\|A\|, \forall A \in \mathbb S^m, P \in {\cal O}^m,
$$
where ${\cal O}^m$ denotes the set of all $m\times m$ orthogonal matrices.
 Corresponding to  $l_p$ norm in $\mathbb R^m$, the Schatten $p$-norm is defined by
$$
\|A\|=\|\lambda (A)\|_p.
$$
Well-known special cases, corresponding to  $l_{\infty}$ norm, the $l_1$ norm and the $l_2$ norm, are  the spectral norm, the nuclear norm  and  the Frobenius norm of a symmetric matrix, respectively.

Many nonlinear optimization problems can be formulated as the following form
\begin{equation}\label{Pconic}
\begin{array}{ll}
\min & f(x) \\[6pt]
{\rm s.t.} & g(x) \in K,
\end{array}
\end{equation}
 where
 $$
 K={\rm epi}\, \|\cdot\|.
 $$
For examples, optimization problems considered in \cite{YLin2007},\cite{BDas2008},\cite{CRecht2009},\cite{CWright2011} and \cite{Qi2013} can be reformulated
 or can be approximated by the mathematical model (\ref{Pconic}), where $\|\cdot\|$ is  spectral or nuclear matrix norm.

Ding, Sun and Toh (2014) \cite{DST14} established several key properties including the closed form solution, calm B-differentiability and strong semi-smoothness of the metric projection operator over the epigraph of  spectral,
and nuclear matrix norm, respectively. Ding (2017) \cite{Ding2017} studied  some variational properties of the Ky Fan $k$-norm $\|\cdot\|_{(k)}$, including  spectral,
and nuclear matrix norm. Ding, Sun, Sun and Toh\cite{DSST2018}\cite{DSST2020} proposed the notion of spectral operator, which is an extension  of L\"{o}wner operator\cite{Lowner1934} of symmetric matrices to non-symmetric matrices and studied variational analysis of spectral operators.

For Problem (\ref{Pconic}), in order to develop optimality conditions, including first-order necessary optimality conditions, second-order necessary optimality conditions and second-order sufficient optimality conditions,  we have to derive the variational geometry of the feasible set $g^{-1}(K)$, namely the tangent cone, the normal cone and the second-order tangent set, this is closely related to the variational geometry of cone $K$. In order to analyze the stability of Problem (\ref{Pconic}), we have to develop the differential properties of projection onto the norm cone,  including directional derivative and B-subdifferential.  How to extend variational properties of the  Ky Fan $k$-norm matrix cone in \cite{Ding2017} to a general norm cone? This motives us to establish the variational analysis for the orthogonally invariant  norm cone of symmetric matrices.

The paper is organized as follows. In Section \ref{Sec:Preliminary}, we give some preliminaries needed in the following sections.
In Section \ref{Section-Normepi}, we develop variational analysis of the orthogonally invariant  norm cone of symmetric matrices, including variational geometry of an orthogonally invariant  norm cone and directional derivative and B-subdifferential of the projection onto the orthogonally invariant  norm cone.
Section \ref{Section-computing} focuses on  calculating $\nabla^2 (\psi\circ \lambda)$ and $G(z)^{-1}$,
which simplifies the formulas for directional derivative and B-subdifferential of the projection operator on the orthogonally invariant  norm cone of symmetric matrices.
In Section \ref{Section-Kn}, we apply the results in Section \ref{Section-Normepi}  and Section \ref{Section-computing}
to Schatten $p$-norm cone when $p \in (1,\infty)$, and the variational geometry of the nuclear norm cone and spectral norm cone is discussed.
We conclude the paper in Section \ref{Section-Conclusion}.
\section{Preliminaries}\label{Sec:Preliminary}
\setcounter{equation}{0}
It follows from \cite[Page 35]{Waston1992} that
all orthogonally invariant  norms of symmetric matrices can be written as
$$
\|A\|=\psi(\lambda (A)),
$$
where $\lambda (A)=(\lambda_1(A),\ldots, \lambda_m(A))^T$ is the vector of eigenvalues of $A \in \mathbb S^m$ with
$$
\lambda_1(A)\geq \lambda_2(A)\geq \cdots\geq  \lambda_m(A),
$$
and $\psi: \mathbb R^m \rightarrow \mathbb R$ is a symmetric gauge function, such a function satisfies the following conditions:
\begin{itemize}
\item[{\rm (i)}] $\psi (y)>0$ for $y \in \mathbb R^m$ with $y\ne 0$;
\item[{\rm (ii)}] $\psi (\alpha y)=|\alpha| \psi (y)$ for $y \in \mathbb R^m$ and $\alpha \in \mathbb R$;
\item[{\rm (iii)}] $\psi ( y^1+y^2)\leq \psi (y^1)+\psi (y^2)$ for $y^1,y^2 \in \mathbb R^m$;
 \item[{\rm (iv)}] $\psi (y_{i_1},\ldots,  y_{i_m})=\psi (y)$ for $y \in \mathbb R^m$ and
     all $i_1,\ldots,i_m$ being a permutation of $1,\ldots, m$.
\end{itemize}
Obviously, $\psi$ is a norm of $\mathbb R^m$. Let $\psi_*$ be the dual norm of $\psi$, namely $\psi_*: \mathbb R^m \rightarrow \mathbb R$,
$$
\psi_*(u)=\displaystyle \max_{\psi(y)\leq 1} \langle u, y\rangle.
$$
 It is easy to check that $\psi_*$ is a symmetric gauge function.

Let $w_1,\ldots, w_r$ be $r$ distinct values of $m$ eigenvalues of $A$,
namely
$$
\lambda_1(A)=\cdots=\lambda_{k_1}(A)=w_1>\lambda_{k_1+1}(A)=\ldots=\lambda_{k_2}(A)=
w_2>\cdots > \lambda_{k_{r-1}+1}=\cdots=\lambda_{k_r}=w_r,
$$
where $k_0=0, k_r=m$. Denote
$$
\alpha_1=\{1,\ldots, k_1\}, \alpha_2=\{k_1+1,\ldots, k_2\}, \ldots, \alpha_r=\{k_{r-1}+1,\ldots, k_r\}.
$$
The following four results are from \cite{LSendov2001}.
\begin{lemma}\label{lem:LS2.1}
Let $\psi: \mathbb R^m \rightarrow \mathbb R$ be a symmetric function, twice differentiable at the point $w \in \mathbb R^m$,
and $P$ be a permutation matrix such that $Pw=w$. Then
\begin{itemize}
\item[{\rm (i)}] $\nabla \psi (w)=P^T\nabla f(w)$ and
\item[{\rm (ii)}]$\nabla^2 \psi(w)=P^T\nabla^2\psi(w)P$.
\end{itemize}
In particular we have the representation
$$
\nabla^2 \psi(w)=\left (
\begin{array}{cccc}
a_{11}\textbf{1}_{|\alpha_1|}\textbf{1}_{|\alpha_1|}^T+b_{k_1}I_{|\alpha_1|} &a_{12}\textbf{1}_{|\alpha_1|}\textbf{1}_{|\alpha_2|}^T&\cdots & a_{1r}\textbf{1}_{|\alpha_1|}\textbf{1}_{|\alpha_r|}^T\\[10pt]
a_{21}\textbf{1}_{|\alpha_2|}\textbf{1}_{|\alpha_1|}^T&a_{22}\textbf{1}_{|\alpha_2|}\textbf{1}_{|\alpha_2|}^T+b_{k_2}I_{|\alpha_2|} &\cdots & a_{2r}\textbf{1}_{|\alpha_2|}\textbf{1}_{|\alpha_r|}^T\\[10pt]
\vdots & \vdots &\vdots& \vdots\\[10pt]
a_{r1}\textbf{1}_{|\alpha_r|}\textbf{1}_{|\alpha_1|}^T&a_{r2}\textbf{1}_{|\alpha_r|}\textbf{1}_{|\alpha_2|}^T & \cdots & a_{rr}\textbf{1}_{|\alpha_r|}\textbf{1}_{|\alpha_r|}^T+b_{k_r}I_{|\alpha_r|}
\end{array}
\right),
$$
where $(a_{ij})_{i,j=1}^r$ is a real symmetric matrix, $b:=(b_1,\ldots, b_m)^T$ is a vector which is block refined by $w$.
\end{lemma}

Let $b(w) \in \mathbb R^m$ with $b(w)=(b_1(w),\ldots, b_m(w))^T$, be defined by
$$
b_i(w)=\left\{
\begin{array}{ll}
\psi''_{ii}(w) & \mbox{if } |I_l|=1,i \in I_l,\\[6pt]
\psi''_{pp}(w)-\psi''_{pq}(w) & \mbox{ for any } p \ne q \in I_l, i \in I_l.
\end{array}
\right.
$$
Define ${\cal A}(w)={\cal A}_{ij}(w)$ by
\begin{equation}\label{eq:calA}
{\cal A}_{ij}(w)=\left\{
\begin{array}{ll}
0 & \mbox{ if } i=j,\\[10pt]
b_i(w) & \mbox{ if } i\ne j \mbox{ but } i, j \in I_l,\\[10pt]
\displaystyle \frac{\psi'_i(w)-\psi'_j(w)}{w_i-w_j} & \mbox{ otherwise}.
\end{array}
\right.
\end{equation}
\begin{lemma}\label{LDlem3.1}
Let $A \in \mathbb S^m$ and suppose $\lambda (A)$ belongs to the domain of the symmetric function $\psi: \mathbb R^m \rightarrow \mathbb R$. Then $\psi$ is differentiable at the point $\lambda(A)$ if and only if $\psi\circ \lambda$ is differentiable at the point $A$. In that case we have the formula
$$
\nabla (\psi \circ \lambda) (A)=U({\rm Diag}\,\nabla \psi(\lambda (A)))U^T
$$
for orthogonal matrix $U$ satisfying $A=U({\rm Diag}\lambda (A))U^T$.
\end{lemma}

\begin{lemma}\label{LDTh3.3}
Let $A \in \mathbb S^m$ and suppose $\lambda (A)$ belongs to the domain of the symmetric function $\psi: \mathbb R^m \rightarrow \mathbb R$. Then $\psi$ is twice differentiable at the point $\lambda(A)$ if and only if $\psi\circ \lambda$ is twice differentiable at the point $A$. Moreover, in this case the Hessian of the spectral function at the matrix $A$ is
$$
\nabla^2 (\psi \circ \lambda) (A)[H]=U\left[{\rm Diag}\,(\nabla^2 \psi(\lambda (A)){\rm diag}\,\widehat H)+{\cal A}\circ \widehat H\right]U^T
$$
for orthogonal matrix $U$ satisfying $A=U({\rm Diag}\lambda (A))U^T$, where $\widehat H=U^THU$ and ${\cal A}={\cal A}(\lambda(A))$ is defined by (\ref{eq:calA}).
\end{lemma}

\begin{lemma}\label{LDTh4.2}
Let $A \in \mathbb S^m$ and suppose $\lambda (A)$ belongs to the domain of the symmetric function $\psi: \mathbb R^m \rightarrow \mathbb R$. Then $\psi$ is twice continuously differentiable at the point $\lambda(A)$ if and only if $\psi\circ \lambda$ is twice  continuously differentiable at the point $A$.
\end{lemma}
\begin{lemma}\label{Fan-Ineq.} (Fan's inequality)\cite{Fan1949}Let $Y,Z \in \mathbb S^m$. Then
$$
\langle Y, Z \rangle\leq
\lambda (Y)^T\lambda (Z),
$$
where the equality holds if and only if $Y$ and $Z$ admit a simultaneous
ordered eigenvalue decomposition, i.e., there exists an orthogonal matrix
$U \in {\cal O}^m$  such that
$$
Y=U{\rm Diag}(\lambda (Y))U^T \mbox{ and } Z=U{\rm Diag}(\lambda (Z))U^T.
$$
\end{lemma}

Let $\psi:\mathbb R^m\rightarrow \mathbb R$ be a symmetric gauge function. The orthogonally invariant norm of real symmetric matrices, denoted by $N$, is defined by
\begin{equation}\label{def:N}
N: \mathbb S^m \rightarrow \mathbb R_+, \,\, N(A)=(\psi\circ \lambda)(A),\, A \in \mathbb S^m.
\end{equation}
\begin{proposition}\label{prop:dual-norm}
The dual norm of $N$, denoted by $N_*$, is expressed as
$$
N_*(A)=(\psi_* \circ \lambda)(A), \quad A \in \mathbb S^m.
$$
\end{proposition}
{\bf Proof.} From the definition of dual norm, we have for $A \in \mathbb S^m$, we have from Lemma \ref{Fan-Ineq.} that
$$
\begin{array}{ll}
N_*(A) & =\displaystyle \max_{N(B) \leq 1}\langle A, B \rangle\\[8pt]
& =\displaystyle \max_{N(B) \leq 1}\langle {\rm Diag}\lambda (A), U^T B U\rangle\\[8pt]
&= \displaystyle \max_{N(U^TBU) \leq 1}\langle {\rm Diag}\lambda (A), U^T B U\rangle\\[8pt]
&= \displaystyle\max_{\psi (y) \leq 1}\langle \lambda (A), y\rangle\\[8pt]
&= \psi_*(\lambda (A))=(\psi_*\circ \lambda)(A).
\end{array}
$$
The proof is completed. \hfill $\Box$\\

The matrix norm cone in $\mathbb S^m$ associated with $N$, denoted by $\K$,  is the epigraph of $N$, namely
\begin{equation}\label{eq:mnormcone}
\K=\{(A,s) \in \mathbb S^{m}\times \mathbb R: N(A) \leq s\}.
\end{equation}
We use  and $\K_*$ to denote the matrix norm cone associated with $N_*$, namely
 the epigraph of $N_*$:
\begin{equation}\label{eq:mdnormcone}
\K_*=\{(A,s) \in  \mathbb S^{m}\times \mathbb R: N_*(A)=(\psi_*\circ \lambda)(A) \leq s\}.
\end{equation}
The polar of ${\cal K}$, denoted by ${\cal K}^{\circ}$, is expressed as
\begin{equation}\label{eq:mpnormcone}
{\cal K}^{\circ}=-{\cal K}_*=\{(A,s) \in  \mathbb S^{m}\times \mathbb R: (\psi_*\circ \lambda)(A)+s \leq 0\}.
\end{equation}
Since a symmetric gauge function is a norm, the following definition is a
variant of differentiable norm  from \cite{Abat1979}.
\begin{definition}\label{sg:psi}
We say that $\psi: \mathbb R^m \rightarrow \mathbb R$ is a differentiable symmetric gauge function if it is differentiable
 at $y \ne 0$. We say $\psi: \mathbb R^m \rightarrow \mathbb R$ is a twice differentiable symmetric gauge function if it is twice differentiable
 at $y \ne 0$.
\end{definition}
In the next section, we will use Morrison formula for linear operators in \cite{DZhou2008}.
\begin{lemma}\label{MRf}(\cite[Theorem 2.1]{DZhou2008})
 Let $A \in B(X), U\in (Y,X)$ and $V\in B(X,Y)$ such that
$A$ is invertible. Then $A+UV$ is invertible if and only if $I_Y +VA^{-1}U$ is invertible. Furthermore,
if $A+UV$ is invertible, then
$$
(A+UV)^{-1}=A^{-1}-A^{-1}U[I_Y +VA^{-1}U]^{-1}VA^{-1}.
$$
\end{lemma}

\section{Variational Analysis of Norm Cone ${\cal K}$}\label{Section-Normepi}
\setcounter{equation}{0}
Let ${\cal H}=\mathbb S^m \times \mathbb R$, define the inner product $\langle \cdot, \cdot \rangle$ in ${\cal H}$:
$$
\langle z, z'\rangle={\rm Tr}(AA')+ss', \, z=(A, s) \in {\cal H},  z'=(A', s') \in {\cal H}.
$$
The induced norm, denoted by $\|\cdot\|_{{\cal H}}$, is defined by
$$
\|z\|_{{\cal H}}=\sqrt{\|A\|_F^2+s^2}.
$$
This section provides basic variational analysis on the convex cone $\K$ mostly through
a study on the orthogonal projection onto $\K$:
\[
  \Pi_{\K} (z) := \argmin_{z' \in \K} \ \displaystyle \frac{1}{2}\| z'-z \|_{{\cal H}}^2.
\]
Let ${\cal I}$ and ${\cal I}'$ be identity mapping in $\mathbb S^m$ and ${\cal H}$, respectively, namely
$$
{\cal I} A=A, \,\, {\cal I}'z=z \mbox{ for } A\in \mathbb S^m \mbox{ and } z \in {\cal H}.
$$

\subsection{The  Tangent Cone and Second-order Tangent Set of $\K$}
 It follows from Lemma \ref{LDlem3.1}, Lemma \ref{LDTh3.3} and Definition \ref{sg:psi} that $N$, defined by (\ref{def:N}), is a  differentiable  norm  if and only if $\psi: \mathbb R^m \rightarrow \mathbb R$  is differentiable
 at $y \ne 0$, and $N$ is a twice differentiable  norm if and only if it is twice differentiable
 at $y \ne 0$.  Obviously, the spectral norm and the nuclear norm  are not differentiable norms, and
 from Lemma \ref{LDlem3.1}, Lemma \ref{LDTh3.3}, the Schatten $p$-norm is differentiable but not twice  differentiable norm
 when $p \in (1,2)$,   the Schatten $p$-norm is a twice differentiable norm when $p \in [2,\infty)$.

\begin{lemma}\label{lem:mtangents}
Let $\psi$  and $\psi_*$ be   differentiable  norms.  Then the following results hold:
\begin{itemize}
\item[{\rm (i)}]  For $(A,s) \in \K$,
$$
{\cal T}_{\K}(A,s)=\left
\{
\begin{array}{ll}
{\cal H} &  N (A)<s,\\[6pt]
\K & (A,s)=(0,0),\\[6pt]
\left\{(d_A,d_s)\in \mathbb{S}^{m} \times \mathbb R: \langle U {\rm Diag}\,(\nabla \psi (\lambda (A)))U^T,d_A \rangle \leq d_s \right\} & N (A)=s>0.
\end{array}
\right.
$$
\item[{\rm (ii)}] For $(A,s) \in \K$,
$$
{\cal N}_{\K}(A,s)=\left
\{
\begin{array}{ll}
\{0\}  &  N (A)<s,\\[6pt]
\K^{\circ} & (A,s)=(0,0),\\[6pt]
\big\{\alpha(Y, -1)\in {\cal H}:(\psi_* \circ \lambda)(Y)=1, \langle Y, A/s\rangle=1,\alpha\geq 0 \big\} & N (A)=s>0.
\end{array}
\right.
$$
\item[{\rm (iii)}] For $(A,s) \in \K^{\circ}$,
$$
{\cal T}_{\K^{\circ}}(A,s)=\left
\{
\begin{array}{ll}
{\cal H} & N_*(A)<-s,\\[6pt]
\K^{\circ} & (A,s)=(0,0),\\[6pt]
\left\{(d_A,d_s)\in {\cal H}: \langle U {\rm Diag}\,(\nabla \psi_* (\lambda (A)))U^T,d_A \rangle \leq -d_s \right\} & N_*(A)=-s>0.
\end{array}
\right.
$$
\item[{\rm (iv)}]
Let $z \in \K$ and $d \in {\cal T}_{ \K}(z)$ where $z=(A,s)$ and $d=(d_A,d_s)$. If $\psi$ is a twice differentiable symmetric gauge function, then
\begin{equation}\label{eq:m2ndTK}
{\cal T}^2_{ \K}(z,d)=\left\{
\begin{array}{ll}
{\cal H} & d \in \mbox{int }{\cal T}_{ \K}(z),\\[10pt]
{\cal T}_{ \K}(z) & z=0,\\[10pt]
\left\{(\xi_A,\xi_s): \begin{array}{l}
\langle U {\rm Diag}\,(\nabla \psi (\lambda (A)))U^T,\xi_A\rangle
-\xi_s\\[6pt]
+\nabla^2 \psi(\lambda (A))[{\rm diag} (\widehat d_A), {\rm diag} (\widehat d_A)]\\[6pt]
+\langle {\cal A}, \widehat d_A \circ \widehat d_A \rangle \leq 0
\end{array}
\right\} & \mbox{otherwise}.
\end{array}
\right.
\end{equation}
\end{itemize}
\end{lemma}
{\bf Proof}. In view of Lemma \ref{LDTh3.3} and Lemma \ref{LDTh4.2}, Assertions (i) and (iv) come from Proposition 2.61 and Proposition 3.30 of \cite{BS00}, respectively. From (\ref{eq:mpnormcone}), and noting Proposition 2.61 of \cite{BS00}, we obtain (iii).  We only need to prove the case when $(\psi \circ \lambda)(A)=s>0$ in (ii).  Since $\K$ is a closed convex cone, we have from (2.110) of \cite{BS00} that
the normal cone of ${\cal K}$ at $z$ is expressed as
\begin{equation}\label{eq:mcone-Normal}
{\cal N}_{\K}(z)=\{(U,t)\in {\cal H}: (U,t)\in \K^{\circ}, \langle (U,t), (A,s)\rangle=0\}.
\end{equation}
Since $\K={\rm epi}\, N={\rm epi}\,(\psi\circ \lambda)$, $(A,s) \in {\rm epi}\, N$, one has from $(U,t)\in {\cal N}_{\K}(z)$ that $t<0$, it is necessary to consider an element $(V,-1) \in {\cal N}_{\K}(z)$. It follows from (\ref{eq:mcone-Normal}) that
\begin{equation}\label{eq:mv-1}
\begin{array}{ll}
(V,-1) \in {\cal N}_{\K}(z)& \Longleftrightarrow (V,-1) \in \K^{\circ}, \langle (V,-1), (A,s)\rangle=0\\[6pt]
&\Longleftrightarrow N_*(V) \leq 1, \langle V, A\rangle=s.
\end{array}
\end{equation}
From the definition of $N_*$, and $N(A/s)=1$ and $\langle V,  A/s\rangle\leq N(A/s)\cdot N_*(V)=N_*(V)$, we have $\langle V,  W/s\rangle \leq 1$. Thus $N_*(V)\leq 1, \langle V, A/s\rangle=1$ is equivalent to $\langle V, W/s\rangle=1$ and $N_*(V)=1$.
The formula of ${\cal N}_{\K}(z)$ for the case  $N(A/s)=1$ follows from (\ref{eq:mv-1}). The proof is completed. \hfill $\Box$
\begin{corollary}\label{cor:NCa}
Let $\psi$  and $\psi_*$ be   differentiable  norms. Then for $(A,s) \in \K$,
$$
{\cal N}_{\K}(A,s)=\left
\{
\begin{array}{ll}
\{0\}  &  N (A)<s,\\[6pt]
\K^{\circ} & (A,s)=(0,0),\\[6pt]
\big\{\alpha(U{\rm Diag}(a)U^T, -1)\in {\cal H}:a\in \partial \psi(\lambda(A)),\alpha\geq 0 \big\} & N (A)=s>0.
\end{array}
\right.
$$
\end{corollary}
{\bf Proof}. We only need to prove the case where $(A,s)$ satisfying $N (A)=s>0$. In this case, we have that $(\psi\circ \lambda)(A)=s$,
 $(\psi_* \circ \lambda)(Y)=1, \langle Y, A/s\rangle=1$,  implying
 $$
 \langle Y, A/s\rangle=N(A/s)N_*(Y).
 $$
 From Lemma \ref{Fan-Ineq.}, we have that  $Y$ and $A$ admit a simultaneous
ordered eigenvalue decomposition, i.e.,
$$
 Y=U{\rm Diag}(\lambda (Y))U^T.
$$
 It follows from Lemma \ref{lem:mtangents} that
 $$
 \begin{array}{ll}
{\cal N}_{\K}(A,s)& =
\big\{\alpha(Y, -1)\in {\cal H}:(\psi_* \circ \lambda)(Y)=1, \langle Y, A/s\rangle=1,\alpha\geq 0 \big\} \\[6pt]
&=\big\{\alpha(U\lambda(Y)U^T, -1):\psi_*(\lambda(Y))=1, \langle \lambda(Y), \lambda(A/s)\rangle=1,\alpha\geq 0 \big\}\\[6pt]
&=\big\{\alpha(U\lambda(Y)U^T, -1):(a,-1) \in {\cal N}_{{\rm epi}\, \psi}(\lambda(A),s) ,\alpha\geq 0 \big\}\\[6pt]
&=\big\{\alpha(U\lambda(Y)U^T, -1):a\in \partial \psi(\lambda(A)), \alpha\geq 0 \big\}.
\end{array}
$$
 The proof is completed. \hfill $\Box$
\subsection{The  Formula for $\Pi_{\K}(z)$}
The projection onto
$\K$ plays an important role in the study of theories and algorithms for the corresponding matrix conic  optimization problems.
In this subsection, we will derive the formula of $\Pi_{\K}(z)$.
For a given  $z=(A,s)\in  {\cal H}$ and any $(w,\mu) \in \mathbb R^{m} \times \mathbb R$ with $w\ne 0$, define $F:\mathbb{R}^{m+1}\rightarrow
\mathbb{R}^{m+1}$ by
\begin{equation}\label{eq:mFp}
F(w,\mu;z)=\left[
\begin{array}{l}
\mu \nabla \psi (w)+w-\lambda (A)\\[3pt]
\psi (w)-\mu -s
\end{array}
\right
].
\end{equation}
\begin{proposition}\label{prop:mproj-f}
Let $\psi: \mathbb R^m \rightarrow \mathbb R_+$ be a differentiable symmetric gauge function. Let $z=(A,s)\in {\cal H}$ be a given point with $A=P{\rm Diag}(\lambda (A))P^T$, then
\begin{equation}\label{eq:mkf}
\Pi_{\K}(A,s)=\left
\{
\begin{array}{ll}
(A,s) & \mbox{if } (A,s) \in \K,\\[2pt]
(0,0) & \mbox{if } (A,s) \in \K^{\circ},\\[2pt]
(P{\rm Diag}(w(z))P^T, \mu(z)+s) &  \mbox{otherwise},
\end{array}
\right.
\end{equation}
where $(w(z),\mu(z))\in \mathbb R^m \times \mathbb R_{++}$ is a solution of the system of equations:
$$
F(w,\mu;z)=0.
$$
\end{proposition}
{\bf Proof}. We may express $\Pi_{{\cal K}}(A,s)$ as the solution of the following optimization problem
\begin{equation}\label{eq:stara}
\min \displaystyle \frac{1}{2} \|X-A\|_F^2+\displaystyle \frac{(t-s)^2}{2}+\delta_{\{(\psi\circ \lambda)(X) \leq t\}}(X,t).
\end{equation}
Noting that
$$
\delta_{\{(\psi\circ \lambda)(X) \leq t\}}(X,t)=\delta_{{\rm epi}\, \psi}(\lambda(X),t)
$$
and from Fy Fan's inequality
$$
\|X-A\|_F \geq \|\lambda (X)-\lambda (A)\|,
$$
we obtain
$$
\begin{array}{l}
\displaystyle \frac{1}{2} \|X-A\|_F^2+\displaystyle \frac{(t-s)^2}{2}+\delta_{\{(\psi\circ \lambda)(X) \leq t\}}(X,t)\\[6pt]
\geq \displaystyle \frac{1}{2} \|\lambda (X)-\lambda (A)\|^2+\displaystyle \frac{(t-s)^2}{2}+\delta_{{\rm epi}\, \psi}(\lambda(X),t).
\end{array}
$$
This implies
$$
\begin{array}{ll}
\displaystyle \frac{1}{2} {\rm dist}^2((A,s),{\cal K})& \geq \displaystyle \frac{1}{2} {\rm dist}^2((\lambda (A),s),{\rm epi}\,\psi)\\[8pt]
&=\displaystyle \frac{1}{2} \|\Pi_{{\rm epi}\,\psi}((\lambda (A),s))-((\lambda (A),s))\|^2.
\end{array}
$$
Let
$$
(\widehat \lambda,\widehat s)=\Pi_{{\rm epi}}(\lambda (A), s),\, \, \widehat X= P{\rm Diag}(\widehat \lambda)P^T.
$$
Then $(\widehat X, \widehat s) \in {\cal K}$, and
$$
\displaystyle \frac{1}{2} \|\widehat X-A\|_F^2+\displaystyle \frac{(\widehat s-s)^2}{2}=\displaystyle \frac{1}{2}
{\rm dist}^2((\lambda (A), s), {\rm epi} \psi).
$$
Therefore, $(\widehat X, \widehat s)$ is the unique solution to Problem
(\ref{eq:stara}).
Now we only need to consider to calculate $\Pi_{{\rm epi}\psi}(\lambda (A), s)$, namely the solution of the following problem
\begin{equation}\label{eq:mP}
\min \displaystyle \frac{1}{2} \|(w,t)-(\lambda (A),s)\|^2 \quad \mbox{s.t. } \psi (w) -t \leq 0.
\end{equation}
It is obvious that $\Pi_{\K}(A,s)=(A,s)$ when $(A,s) \in \K$ and $\Pi_{\K}(A,s)=(0,0)$
when $(A,s) \in [\K]^{\circ}$. If $(A,s) \notin \left( \K\cup  [\K]^{\circ}\right)$, then for
$(A_*,s_*)=\Pi_{\K}(A,s)$, one has that $A_* \ne 0$ and  $s_*>0$. Obviously, one has that $(A_*,s_*) \in {\rm bdry}\,\K$, namely $N(A_*)-s_*=0$. From the above analysis, we know that $\Pi_{\K}(A,s)=(P{\rm Diag}(w^*)P^T,s^*)$ where $(w_*,s_*)$ is the unique solution to the  convex programming problem (\ref{eq:mP}).

Let ${\cal L}: \mathbb R^{m+1} \times \mathbb R \rightarrow \mathbb R$ be the Lagrangian of Problem (\ref{eq:mP})
$$
{\cal L}(w,t,\mu;z)=\displaystyle \frac{1}{2} \|(w,t)-(\lambda(A),s)\|^2+\mu(\psi(w) -t ).
$$
Then $(w_*,s_*)$ is a solution to  Problem (\ref{eq:mP}) if and only if there exists a real number $\mu_* \geq 0$ such that
$(w_*,s_*,\mu_*)$ satisfies
$$
\begin{array}{l}
\nabla_w {\cal L}(w,t,\mu;z)=0,\\[2pt]
\displaystyle\frac{\partial {\cal L}}{\partial t}(w,t,\mu;z)=0,\\[2pt]
\mu(\psi(w)-t)=0, \mu \geq 0,
\end{array}
$$
which are equivalent to
$$
\begin{array}{l}
\mu \nabla \psi (w)+w-\lambda (A)=0,\\[2pt]
t-s-\mu=0,\\[2pt]
\mu (\psi(w)-t)=0, \mu \geq 0.
\end{array}
$$
From the above system, we know $\mu_* >0$, because otherwise $\Pi_{\K}(A,s)=(A,s)$, which contradicts with the assumption $(A,s) \notin \left( \K\cup  [\K]^{\circ}\right)$. Therefore, $(w_*,s_*,\mu_*)$ satisfies
$$
\begin{array}{l}
\mu \nabla \psi(w)+w-\lambda (A)=0,\\[2pt]
t-s-\mu=0,\\[2pt]
\psi (w)-s-\mu=0, \mu >0.
\end{array}
$$
Namely $(w_*,s_*)=(w(z),\mu(z)+s)$, where $(w(z),\mu(z))\in \mathbb R^m \times \mathbb R_{++}$ is a solution of the system of equations $
F(w,\mu;z)=0$. The proof is completed. \hfill $\Box$\\
We use $\K_2$ to denote the second-order cone in $\mathbb S^m \times \mathbb R$, namely the matrix cone associated with
$$\psi (w)=\left(\displaystyle \sum_{j=1}^m w_j^2\right)^{1/2}=\|w\|_2,
$$
\begin{equation}\label{eq:m2ndcone}
\K_2=\{(A,s) \in \mathbb S^m \times \mathbb R: \|A\|_F \leq s\}.
\end{equation}
Then we obtain, from Proposition \ref{prop:mproj-f}, the formula of the projection operator onto $\K_2$.
\begin{corollary}\label{prop:mproj-fcor}
Let $z=(A,s)\in {\cal H}$ be a given point, then
\begin{equation}\label{eq:mkf2}
\Pi_{\K_2}(A,s)=\left
\{
\begin{array}{ll}
(A,s) & \mbox{if } (A,s) \in \K_2,\\[2pt]
(0,0) & \mbox{if } (A,s) \in [\K_2]^{\circ},\\[2pt]
(\|A\|_F+s)\cdot\displaystyle \frac{1}{2}\left(
\displaystyle \frac{A}{\|A\|_F},1\right) &  \mbox{otherwise}.
\end{array}
\right.
\end{equation}
\end{corollary}
{\bf Proof}. It follows from Proposition \ref{prop:mproj-f}, we only need to solve the system of equations
$F(w,\mu;z)=0$.  In this case, this system becomes
$$
\left
\{
\begin{array}{l}
\mu \cdot \displaystyle \frac{w}{\|w\|_2}+w-\lambda (A)=0,\\[12pt]
\|w\|_2-\mu-s=0.
\end{array}
\right.
$$
Solving this system of equations, we obtain
$$
w=\displaystyle \frac{1}{2}\cdot (\|\lambda (A)\|_2+s)\cdot
\displaystyle \frac{\lambda (A)}{\|\lambda (A)\|_2},\quad \mu=\displaystyle \frac{1}{2}\cdot (\|\lambda (A)\|_2-s),
$$
namely
$$
w=\displaystyle \frac{1}{2}\cdot (\|A\|_F+s)\cdot
\displaystyle \frac{\lambda (A)}{\|A)\|_F},\quad \mu=\displaystyle \frac{1}{2}\cdot (\|A\|_F-s).
$$
Therefore, we get
$$
\begin{array}{l}
(P{\rm Diag}(w)P^T, \mu+s)\\[10pt]
=\left(\displaystyle \frac{1}{2}\cdot (\|A\|_F+s)\cdot
\displaystyle \frac{P\lambda (A)P^T}{\|A)\|_F},\displaystyle \frac{1}{2}\cdot (\|A\|_F+s)\right)\\[10pt]
=(\|A\|_F+s)\cdot\displaystyle \frac{1}{2}\left(
\displaystyle \frac{A}{\|A\|_F},1\right),
\end{array}
$$
which yields the result desired. \hfill $\Box$

\subsection{Directional Derivative and B-subdifferential of $\Pi_{\K}(z)$}
Since $\Pi_{\K}(z)$ is globally Lipschitz continuous with constant $1$, it is differentiable almost everywhere. From the expression of $\Pi_{\K}(z)$, we know that it is differentiable at points in the following set
$$
{\rm int}\, \K\cup {\rm int}\,[\K]^{\circ} \cup {\rm int}\,[{\cal H}\setminus (\K\cup \K^{\circ})].
$$
Now we give the formulas for the derivatives of  $\Pi_{\K}$ at points in the above set.
\begin{proposition}\label{mpropPiD}
Let $\psi: \mathbb R^m \rightarrow \mathbb R_+$ be a twice differentiable symmetric gauge function. Let $z=(A,s)\in {\cal H}$ be a given point with $A=P{\rm Diag}(\lambda (A))P^T$, then
  the following results hold:
\begin{itemize}
\item[{\rm (i)}] If $z=(A,s) \in {\rm int } \K$, then ${\rm D}\Pi_{\K}(z)(d_A,d_s)=(d_A,d_s)$ for
$(d_A,d_s) \in {\cal H}$;
\item[{\rm (ii)}] If $z=(A,s) \in {\rm int } \K^{\circ}$, then ${\rm D}\Pi_{\K}(z)=(d_A,d_s)=0$ for
$(d_A,d_s) \in {\cal H}$;
\item[{\rm (iii)}] If $z=(A,s) \in {\rm int }[{\cal H}\setminus (\K\cup[\K]^{\circ})]$, then for
$(d_A,d_s) \in {\cal H}$,
\begin{equation}\label{eq:mprojD}
\begin{array}{l}
{\rm D}\Pi_{\K}(z)\left[
\begin{array}{l}
\delta_A\\[8pt]
\delta_s
\end{array}
\right] =
\left[
\begin{array}{ll}
G(z)^{-1} & 0\\[8pt]
0 & 1
\end{array}
\right] \left[
\begin{array}{l}
\delta_A\\[8pt]
\delta_s
\end{array}
\right] \\[4pt]
\quad -\beta(z)\left[
\begin{array}{c}
G(z)^{-1}\nabla (\psi\circ \lambda)(X(z))\\[8pt]
-1
\end{array}
\right][\nabla (\psi\circ \lambda)(X(z))^*G(z)^{-1} \,\,\, -1]\left[
\begin{array}{l}
\delta_A\\[8pt]
\delta_s
\end{array}
\right],
\end{array}
\end{equation}
where  $$\beta(z)=(1+\nabla (\psi\circ \lambda)(X(z))^*G(z)^{-1}\nabla (\psi\circ \lambda)(X(z)))^{-1},\, G(z)={\cal I}+\mu(z) \nabla^2 (\psi\circ \lambda)(X(z)),$$
where  $X(z)=P {\rm Diag}(u(z))P^T$ and $(u(z),\mu(z))\in \mathbb R^m \times \mathbb R_{++}$  is a solution of the system of equations:
$$
F(u,\mu;z)=0.
$$
\end{itemize}
\end{proposition}
{\bf Proof}. From the expression of $\Pi_{\K}(z)$, the results for (i) and (ii) are obvious, we only need to verify the result for the case (iii).

Let $(X(z),t(z)) =\Pi_{\K}(z)$, then $(X(z),t(z)) =(P{\rm Diag}(u(z))P^T,t(z))$ with $(u(z),t(z))=\Pi_{{\rm epi}\,\psi}((\lambda (A),s))$
and $((u(z),t(z)),\mu(z))$ is the Karush-Kuhn-Tucker pair of the following problem
\begin{equation}\label{eq:mPa}
\min \displaystyle \frac{1}{2} \|(u,t)-(\lambda (A),s)\|^2 \quad \mbox{s.t. } \psi (u) -t \leq 0
\end{equation}
if and only if $((X(z),t(z)),\mu(z))$ is the Karush-Kuhn-Tucker pair of the following problem
\begin{equation}\label{eq:starab}
\begin{array}{ll}
\min  & \displaystyle \frac{1}{2} \|X-A\|_F^2+\displaystyle \frac{(t-s)^2}{2}\\[6pt]
{\rm s.t.} & (\psi\circ \lambda)(X) \leq t.
\end{array}
\end{equation}
It follows from Proposition \ref{prop:mproj-f} that $(P{\rm Diag}(u(z))P^T,t(z),\mu(z))$ is  the Karush-Kuhn-Tucker pair
of Problem (\ref{eq:starab}) if and only if $X(z)=P{\rm Diag}(u(z))P^T$,$t(z)=s+\mu(z)$ and $(u(z),\mu(z))$ is the solution of the following system
$$
F(u,\mu;z)=0.
$$
Let the Lagrangian of Problem (\ref{eq:starab}) be $L(X,t,\mu)$,
$$
L(X,t,\mu)=\displaystyle \frac{1}{2} \|X-A\|_F^2+\displaystyle \frac{(t-s)^2}{2}+\mu((\psi \circ\lambda)(X)-t).
$$
The optimality conditions of Problem (\ref{eq:starab}) can be expressed as
$$
\begin{array}{l}
{\rm D}_X L(X,t,\mu)=0,\\[6pt]
\partial_t  L(X,t,\mu)=0,\\[6pt]
0\leq \mu, (\psi\circ{\lambda})(X)-t \leq 0, \mu[(\psi\circ{\lambda})(X)-t]=0.
\end{array}
$$
Since $z=(A,s) \in {\rm int }[{\cal H}\setminus (\K\cup[\K]^{\circ})]$, we can easily check
$\mu(z)>0$, this yields
$$
\begin{array}{l}
X-A+\mu \nabla (\psi \circ \lambda) (X)=0,\\[6pt]
t-s-\mu=0,\\[6pt]
(\psi\circ{\lambda})(X)-t=0.
\end{array}
$$
or
\begin{equation}\label{eq:Optimality-kkt}
{\cal F}(X,\mu;z)=\left(
\begin{array}{l}
X-A+\mu \nabla (\psi \circ \lambda) (X)\\[6pt]
(\psi\circ{\lambda})(X)-s-\mu
\end{array}
\right)=0.
\end{equation}
It is easy to check
$$
D_{X,\mu}{\cal F}(X,\mu;z)=
\left[
\begin{array}{cc}
{\cal I} +\mu \nabla^2( \psi\circ \lambda) (X) & \nabla (\psi\circ \lambda)(X)\\[3pt]
\nabla (\psi\circ \lambda)(X)^* & -1
\end{array}
\right].
$$
Since the Schur complement matrix
$$
D_{X,\mu}{\cal F}(X,\mu;z)/-1={\cal I} +\mu \nabla^2( \psi\circ \lambda) (X)+ \nabla (\psi\circ \lambda)(X) \nabla (\psi\circ \lambda)(X)^*
$$
is self-adjoint and positively definite, the operator $D_{X,\mu}{\cal F}(X,\mu;z)$ is nonsingular.  It follows from
$z=(A,s) \in {\rm int }[\mathbb R^{m+1}\setminus (\K\cup[\K]^{\circ})]$ and the classical implicit function theorem that there are exists $\delta>0$ and $\varepsilon>0$ and a mapping $(\tilde X(\cdot),\tilde \mu(\cdot)):\textbf{B}_{\delta}(z)\rightarrow
\textbf{B}_{\varepsilon}(X(z),\mu(z))$ such that $\tilde X(z)=X(z)$ and $\tilde \mu(z)=\mu(z)$, and for any $z'=(A',s') \in \textbf{B}_{\delta}(z)$, $(\tilde X(z'),\tilde t(z'),\tilde \mu(z'))$ satisfies ${\cal F}(\tilde X(z'),\tilde \mu(z');z')=0$ and  $\tilde \mu(z')=\tilde t(z')-s'$  and they are
continuously differentiable at $z'=z$. Taking derivative on the both sides of ${\cal F}(\tilde X(z'),\tilde \mu(z');z')=0$ with respect to $z'$ at $z'=z$ along $\delta z=(\delta_A,\delta_s)$, we obtain for $\tilde X':={\rm D}\tilde X(z)\delta z$ and
$\tilde \mu':={\rm D}\tilde \mu(z)\delta z$ that
$$
\begin{array}{ll}
X'-\delta_A+\mu(z)\nabla^2 (\psi\circ \lambda)(X(z))[X']+\mu' \nabla (\psi\circ \lambda)(X(z))=0,\\[8pt]
\langle \nabla (\psi\circ \lambda)(X(z)), X'\rangle  -\delta_s-\mu'=0.
\end{array}
$$
The above equations can be written as
$$
\left[
\begin{array}{cc}
{\cal I}+\mu(z) \nabla^2 (\psi\circ \lambda)(X(z)) & \nabla (\psi\circ \lambda)(X(z))\\[6pt]
\nabla (\psi\circ \lambda)(X(z))^* & -1
\end{array}
\right]\left[
\begin{array}{l}
X'\\[6pt]
\mu'
\end{array}
\right]=\left[
\begin{array}{l}
\delta_A\\[6pt]
\delta_s
\end{array}
\right],
$$
or  equivalently
$$
D_{X,\mu}{\cal F}(X(z),\mu(z);z)\left[
\begin{array}{l}
X'\\[6pt]
\mu'
\end{array}
\right]=\left[
\begin{array}{l}
\delta_A\\[6pt]
\delta_s
\end{array}
\right].
$$
Since $G(z)$ is self-adjoint and positively definite and  the Schur complement matrix
$$
D_{X,\mu}{\cal F}(X(z),\mu(z);z)/-1=G(z)+ \nabla (\psi\circ \lambda)(X(z)) \nabla (\psi\circ \lambda)(X(z))^*
$$
is self-adjoint and positively definite, the operator $D_{X,\mu}{\cal F}(X(z),\mu(z);z)$ is nonsingular. Thus we obtain
$$
\left[
\begin{array}{l}
X'\\[6pt]
\mu'
\end{array}
\right]=\left[
\begin{array}{cc}
G(z) & \nabla (\psi\circ \lambda)(X(z))\\[6pt]
\nabla (\psi\circ \lambda)(X(z))^* & -1
\end{array}
\right]^{-1}\left[
\begin{array}{l}
\delta_A\\[6pt]
\delta_s
\end{array}
\right].
$$
From the Sherman-Morrison-Woodbury formula of linear operators between Banach spaces in \cite{DZhou2008}, we have
$$
\begin{array}{l}
[G(z)+ \nabla (\psi\circ \lambda)(X(z)) \nabla (\psi\circ \lambda)(X(z))^*]^{-1}\\[8pt]
=G(z)^{-1}-\beta (z)G(z)^{-1}\nabla (\psi\circ\lambda)(X(z))\nabla (\psi\circ\lambda)(X(z))^*G(z)^{-1},
\end{array}
$$
where $(1+\nabla (\psi\circ \lambda)(X(z))^*G(z)^{-1}\nabla (\psi\circ \lambda)(X(z)))^{-1}$.
Therefore we obtain
$$
\begin{array}{l}
\left[
\begin{array}{l}
X'\\[6pt]
\mu'
\end{array}
\right]=\left[
\begin{array}{cc}
G(z) & \nabla (\psi\circ \lambda)(X(z))\\[6pt]
\nabla (\psi\circ \lambda)(X(z))^* & -1
\end{array}
\right]^{-1}\left[
\begin{array}{l}
\delta_A\\[6pt]
\delta_s
\end{array}
\right]\\[20pt]
=\left[
\begin{array}{ll}
[G(z)+ \nabla (\psi\circ \lambda)(X(z)) \nabla (\psi\circ \lambda)(X(z))^*]^{-1}& \beta (z)G(z)^{-1}\nabla (\psi\circ\lambda)(X(z))\\[14pt]
\beta (z)\nabla (\psi\circ\lambda)(X(z))^*G(z)^{-1} &-\beta (z)
\end{array}
\right]\left[
\begin{array}{l}
\delta_A\\[6pt]
\delta_s
\end{array}
\right].
\end{array}
$$
Noting that $t(z)=\mu(z)+s$, we obtain
$$
\begin{array}{l}
{\rm D} \Pi_{\K}(z)\delta_z=\left[
\begin{array}{l}
X'\\[3pt]
\mu'
\end{array}
\right]+\left[
\begin{array}{l}
0\\[3pt]
\delta_s
\end{array}
\right]\\[16pt]
=\left[
\begin{array}{c}
G(z)^{-1}d_A\\[3pt]
\delta_s
\end{array}
\right]-\beta(z)\left[
\begin{array}{c}
G(z)^{-1}\nabla (\psi\circ \lambda)(X(z))\\[3pt]
-1
\end{array}
\right](\nabla (\psi\circ \lambda)(X(z))^*G(z)^{-1}\delta_A-\delta_s)\\[16pt]
=\left[
\begin{array}{ll}
G(z)^{-1} & 0\\[14pt]
0 & 1
\end{array}
\right]\left[
\begin{array}{c}
d_A\\[3pt]
\delta_s
\end{array}
\right] \\[4pt]
\quad -\beta(z)\left[
\begin{array}{c}
G(z)^{-1}\nabla (\psi\circ \lambda)(X(z))\\[14pt]
-1
\end{array}
\right][\nabla (\psi\circ \lambda)(X(z))^*G(z)^{-1} \,\,\, -1]\left[
\begin{array}{c}
d_A\\[3pt]
\delta_s
\end{array}
\right],
\end{array}
$$
which yields (\ref{eq:mprojD}). \hfill $\Box$\\
Define
$$
F_*(u,\mu;z)=\left[
\begin{array}{l}
\mu \nabla \psi_* (u)+u-\lambda (A)\\[3pt]
\psi_*(u)-\mu +s
\end{array}
\right
].
$$
From the expression (\ref{eq:mpnormcone}) of $\K^{\circ}$, like Proposition \ref{mpropPiD}, we can easily obtain the following result.
\begin{proposition}\label{mpropqPiDq}
Let $\psi_*: \mathbb R^m \rightarrow \mathbb R_+$ be a twice differentiable symmetric gauge function. Let $z=(A,s)\in {\cal H}$ be a given point with $A=P{\rm Diag}(\lambda (A))P^T$, then
 the following results hold:
\begin{itemize}
\item[{\rm (i)}] If $z=(A,s) \in {\rm int } \K^{\circ}$, then ${\rm D}\Pi_{\K^{\circ}}(z)(d_A,d_s)=(d_A,d_s)$ for
$(d_A,d_s) \in {\cal H}$;
\item[{\rm (ii)}] If $z=(A,s) \in {\rm int } \K$, then ${\rm D}\Pi_{\K^{\circ}}(z)=(d_A,d_s)=0$ for
$(d_A,d_s) \in {\cal H}$;
\item[{\rm (iii)}] If $z=(A,s) \in {\rm int }[{\cal H}\setminus (\K\cup[\K]^{\circ})]$, then for
$(d_A,d_s) \in {\cal H}$,
\begin{equation}\label{meq:projD}
\begin{array}{l}
{\rm D}\Pi_{\K^{\circ}}(z)\left[
\begin{array}{l}
\delta_A\\[8pt]
\delta_s
\end{array}
\right] =
\left[
\begin{array}{ll}
G_*(z)^{-1} & 0\\[8pt]
0 & 1
\end{array}
\right] \left[
\begin{array}{l}
\delta_A\\[8pt]
\delta_s
\end{array}
\right] \\[4pt]
\quad -\beta_*(z)\left[
\begin{array}{c}
G_*(z)^{-1}\nabla (\psi_*\circ \lambda)(X(z))\\[8pt]
1
\end{array}
\right][\nabla (\psi_*\circ \lambda)(X(z))^*G_*(z)^{-1} \,\,\, 1]\left[
\begin{array}{l}
\delta_A\\[8pt]
\delta_s
\end{array}
\right],
\end{array}
\end{equation}
where  $$\beta_*(z)=(1+\nabla (\psi_*\circ \lambda)(X(z))^*G_*(z)^{-1}\nabla (\psi_*\circ \lambda)(X(z)))^{-1},\, G_*(z)={\cal I}+\mu(z) \nabla^2 (\psi_*\circ \lambda)(X(z)),$$
where ${\cal I}: \mathbb S^m \rightarrow \mathbb S^m$ is the identity mapping, $X(z)=P {\rm Diag}(u(z))P^T$ and $(u(z),\mu(z))\in \mathbb R^m \times \mathbb R_{++}$ being a solution of the system of equations:
$$
F_*(u,\mu;z)=0.
$$
\end{itemize}
\end{proposition}

For the second-order cone in $\mathbb S^{m}\times \mathbb R$, we have the following results about the derivatives of the projection operator
on to $\K_2$.
\begin{corollary}\label{coromPiD2}
Let $z=(A,s)\in {\cal H}$ be a given point with $A=P{\rm Diag}(\lambda (A))P^T$, then
 Then the following results hold:
\begin{itemize}
\item[{\rm (i)}] If $z=(A,s) \in {\rm int } \K_2$, then ${\rm D}\Pi_{\K_2}(z)(d_A,d_s)=(d_A,d_s)$ for
$(d_A,d_s) \in {\cal H}$;
\item[{\rm (ii)}] If $z=(A,s) \in {\rm int } \K_2^{\circ}$, then ${\rm D}\Pi_{\K_2}(z)=(d_A,d_s)=0$ for
$(d_A,d_s) \in {\cal H}$;
\item[{\rm (iii)}] If $z=(A,s) \in {\rm int }[{\cal H}\setminus (\K_2\cup[\K_2]^{\circ})]$, then
for
$(d_A,d_s) \in {\cal H}$,
\begin{equation}\label{meq:projD2}
{\rm D}\Pi_{\K_2}(z)\left[
\begin{array}{l}
\delta_A\\[20pt]
\delta_s
\end{array}
\right]=\displaystyle \frac{1}{2}\left[
\begin{array}{rc}
{\cal I}+\displaystyle \frac{s}{\|A\|_F}{\cal I}-\displaystyle \frac{s}{\|A\|_F}\displaystyle \frac{AA^*}{\|A\|_F^2} & \displaystyle \frac{A}{\|A\|_F}\\[14pt]
\displaystyle \frac{A^*}{\|A\|_F} & 1
\end{array}
\right]\left[
\begin{array}{l}
\delta_A\\[20pt]
\delta_s
\end{array}
\right].
\end{equation}
\end{itemize}
\end{corollary}
{\bf Proof}. The results in (i) and (ii) are obvious. We only need to prove (iii).
As in the proof of Corollary \ref{prop:mproj-fcor}, it is easy to obtain
$$
X(z)=\displaystyle \frac{1}{2}\cdot (\|A\|_F+s)\cdot
\displaystyle \frac{A}{\|A\|_F},\quad \mu(z)=\displaystyle \frac{1}{2}\cdot (\|A\|_F-s)
$$
for $z=(A,s)$.  Noting for $X(z) \ne 0$, $\psi$ is twice continuously differentiable with
$$
\nabla [\psi \circ \lambda](X(z))=\displaystyle \frac{X(z)}{\|X(z)\|_F}=\displaystyle \frac{A}{\|A\|_F}.
$$
For ${\cal A}$ defined by (\ref{eq:calA}), it is easy to verify
\begin{equation}\label{eq:l2calA}
{\cal A}=\textbf{1}_m\textbf{1}_m^T-I_m.
\end{equation}
And for $u(z)$, we have
\begin{equation}\label{eq:uTp}
u(z)^T
\left
[
\begin{array}{c}
p_1^THp_1\\[3pt]
\vdots\\[3pt]
p_m^THp_m
\end{array}
\right]=\displaystyle \sum_{i=1}^m u_i(z)p_i^THp_i=\displaystyle \sum_{i=1}^m\langle u_i(z)p_ip_i^T,H\rangle=\langle X(z), H\rangle.
\end{equation}
Combing (\ref{eq:l2calA}) and (\ref{eq:uTp}), we have from Lemma \ref{LDTh3.3} that
$$
\begin{array}{l}
\nabla^2 (\psi \circ \lambda)(X(z))[H]\\[6pt]
=P\left[ {\rm Diag}\left( \displaystyle \frac{1}{\|u(z)\|_2}\left(I-\displaystyle \frac{u(z)u(z)^T}{\|u(z)\|_2^2}\right)\left
[
\begin{array}{c}
p_1^THp_1\\[3pt]
\vdots\\[3pt]
p_m^THp_m
\end{array}
\right]\right)+{\cal A} \circ \widehat H\right]P^T\\[10pt]
=P \left[ \displaystyle \frac{1}{\|u(z)\|_2}I \circ \widehat H -\displaystyle \frac{\langle X(z), H\rangle}{\|u(z)\|^3}
{\rm Diag} (u)+{\cal A}\circ \widehat H\right] P^T\\[10pt]
=P \left[ \displaystyle \frac{1}{\|u(z)\|_2} \widehat H -\displaystyle \frac{\langle X(z), H\rangle}{\|u(z)\|^3}
{\rm Diag} (u(z))\right] P^T\\[10pt]
=\displaystyle \frac{1}{\|u(z)\|_2}  H-\displaystyle \frac{\langle X(z), H\rangle}{\|u(z)\|^3}X(z)\\[10pt]
=\displaystyle \frac{1}{\|u(z)\|_2} \left[ {\cal I}-\displaystyle \frac{X(z)X(z)^*}{\|u(z)\|^2}\right] H.
\end{array}
$$
We obtain
$$
\begin{array}{ll}
G(z) &= {\cal I}+ \mu (z) \nabla^2(\psi\circ \lambda)(X(z))\\[8pt]
&={\cal I} +\displaystyle \frac{\mu(z)}{\|u(z)\|_2} \left[ {\cal I}-\displaystyle \frac{X(z)X(z)^*}{\|u(z)\|^2}\right]\\[12pt]
&= \displaystyle \frac{2\|A\|_F}{\|A\|_F+s} \left({\cal I}- \displaystyle \frac{\|A\|_F-s}{2\|A\|_F}\cdot  \displaystyle \frac{AA^*}{\|A\|_F^2}\right).
\end{array}
$$
By using Sherman-Morrison formula, we obtain
$$
G(z)^{-1}=\displaystyle \frac{1}{2}\left(1+ \displaystyle \frac{s}{\|A\|_F} \right){\cal I}+\displaystyle \frac{1}{2}\left(1- \displaystyle \frac{s}{\|A\|_F} \right)\displaystyle \frac{AA^*}{\|A\|_F^2}.
$$
Therefore, we obtain
$$
G(z)^{-1} \nabla (\psi\circ \lambda)(X(z))=G(z)^{-1}\displaystyle \frac{A}{\|A\|_F}=\displaystyle \frac{A}{\|A\|_F}.
$$
and
$$
 \nabla (\psi\circ \lambda)(X(z))^*G(z)^{-1} \nabla (\psi\circ \lambda)(X(z))=1,
$$
implying
$$
(1+\nabla (\psi\circ \lambda)(X(z))^*G(z)^{-1} \nabla (\psi\circ \lambda)(X(z)))^{-1}=1/2.
$$
Thus we get
$$
\begin{array}{l}
G(z)^{-1}-\beta(z)G(z)^{-1}\nabla (\psi\circ \lambda)(X(z))\nabla (\psi\circ \lambda)(X(z))^*G(z)^{-1}\\[10pt]
=\displaystyle \frac{1}{2}\left[{\cal I}+\displaystyle \frac{s}{\|A\|_F}{\cal I}
-\displaystyle \frac{s}{\|A\|_F}\displaystyle \frac{AA^*}{\|A\|_F^2}\right].
 \end{array}
$$
Therefore we obtain
$$
\begin{array}{ll}
{\rm D} \Pi_{{\cal K}}(z) &=\left[  \begin{array}{ll}
\displaystyle \frac{1}{2}\left(1+ \displaystyle \frac{s}{\|A\|_F} \right){\cal I}+\displaystyle \frac{1}{2}\left(1- \displaystyle \frac{s}{\|A\|_F} \right)\displaystyle \frac{AA^*}{\|A\|_F^2} & 0\\[10pt]
0 & 1
\end{array}
\right]-\displaystyle \frac{1}{2}\left[
\begin{array}{c}
\displaystyle \frac{A}{\|A\|_F}\\[6pt]
-1
\end{array}
\right]\left [ \displaystyle \frac{A^*}{\|A\|_F} \,\, -1\right]\\[16pt]
&=\displaystyle \frac{1}{2}\left[
\begin{array}{rc}
{\cal I}+\displaystyle \frac{s}{\|A\|_F}{\cal I}-\displaystyle \frac{s}{\|A\|_F}\displaystyle \frac{AA^*}{\|A\|_F^2} & \displaystyle \frac{A}{\|A\|_F}\\[14pt]
\displaystyle \frac{A^*}{\|A\|_F} & 1
\end{array}
\right].
\end{array}
$$
Namely (\ref{meq:projD2}) holds.  The proof is completed. \hfill $\Box$\\
In the following theorem, we will derive the formulas of directional derivatives of $\Pi_{\K}$, which  are very important in stability analysis for Problem (\ref{Pconic}).
\begin{theorem}\label{mthPiDirction}
Let $\psi$ be a  differentiable symmetric gauge function. Let $z=(A,s)\in {\cal H}$ be a given point with $A=P{\rm Diag}(\lambda (A))P^T$, then the following results hold:
\begin{itemize}
\item[{\rm (i)}] If $z=(A,s) \in {\rm int } \K$, then $\Pi'_{\K}(z;d_z)=d_z$;
\item[{\rm (ii)}] If $z=(A,s) \in {\rm int } [\K]^{\circ}$, then $\Pi'_{\K}(z;d_z)=0$;
\item[{\rm (iii)}] If $z=(A,s) \in {\rm int }[\mathbb S^{m}\times \mathbb R\setminus (\K\cup[\K]^{\circ})]$ and
$\psi$ is a twice  differentiable   symmetric gauge function, then
\begin{equation}\label{meq:projD-1}
\begin{array}{l}
\Pi'_{\K}(z;d_z)=
\left[
\begin{array}{c}
G(z)^{-1}d_A\\[14pt]
d_s
\end{array}
\right] \\[4pt]
\quad -\displaystyle \frac{[\nabla (\psi \circ \lambda)(X(z))^*G(z)^{-1}d_A -d_s]}{1+\nabla (\psi \circ \lambda)(X(z))^*G(z)^{-1}\nabla (\psi \circ \lambda)(X(z))}\left[
\begin{array}{c}
G(z)^{-1}\nabla (\psi \circ \lambda)(X(z))\\[14pt]
-1
\end{array}
\right],
\end{array}
\end{equation}
where $X(z)=P{\rm Diag}u(z)P^T$, $G(z)={\cal I}+\mu(z) \nabla^2 (\psi\circ\lambda)(X(z))$ with $(u(z),\mu(z))\in \mathbb R^m \times \mathbb R_{++}$ being a solution of the system of equations:
$$
F(u,\mu;z)=0.
$$
\item[{\rm (iv)}] If $z=(A,s) \in {\rm bdry} \K \setminus\{0\}$ and
$\psi$ is a twice  differentiable   symmetric gauge function, then
  \begin{equation}\label{meq:projD-2}
\begin{array}{l}
\Pi'_{\K}(z;d_z)=
\left[
\begin{array}{c}
d_A\\[14pt]
d_s
\end{array}
\right]-\displaystyle \frac{[\nabla (\psi\circ \lambda)(A)^*d_A -d_s]_+}{1+\|\nabla (\psi\circ \lambda)(A)\|^2_F}\left[
\begin{array}{c}
\nabla (\psi\circ \lambda)(A)\\[14pt]
-1
\end{array}
\right].
\end{array}
\end{equation}
\item[{\rm (v)}] If $z=(A,s) \in {\rm bdry} [\K]^{\circ} \setminus\{0\}$ and
$\psi_*$ be  twice differentiable symmetric gauge function, then
  \begin{equation}\label{eq:projD-3}
\begin{array}{l}
\Pi'_{\K}(z;d_z)=
\displaystyle \frac{[\nabla (\psi_*\circ \lambda)(A)^*d_A +d_s]_+}{1+\|\nabla (\psi_*\circ \lambda)(A)\|^2_F}\left[
\begin{array}{c}
\nabla (\psi_*\circ \lambda)(A)\\[14pt]
1
\end{array}
\right].
\end{array}
\end{equation}
\item[{\rm (vi)}] If $z=(A,s)=(0,0)$, then $\Pi'_{\K}(z;d_z)=\Pi_{\K}({\rm d}_z)$.
\end{itemize}
\end{theorem}
{\bf Proof}. Assertions (i),(ii),(iii),(vi) are obvious. We only need to (iv) and (v). For case (iv), $z=\Pi_{{\cal K}}(z)$, if $d_z \in {\cal R}_{\K}(z)$,\footnote{${\cal R}_C(a)$ denotes the radial cone of a nonempty convex set $C$ at $a\in C$, which defined by
$${\cal R}_C(a)=\bigcup_{\lambda \geq 0, a'\in C}\{\lambda (a'-a)\}.$$
}
then $\nabla (\psi \circ \lambda)(A)^*d_A -d_s\leq 0$, $\Pi'_{\K}(z;d_z)=d_z$ and (\ref{meq:projD-2}) is satisfied for this case.  If $d_z \in {\cal T}_{\K}(z)$, then there exists a sequence $d^k \in {\cal R}_{\K}(z)$ such that
$d^k \rightarrow d_z$. Noting $\Pi'_{\K}(z;d^k)=d^k$ and $d \rightarrow \Pi'_{\K}(z;d)$ is Lipschitz continuous, we obtain $\Pi'_{\K}(z;d_z)=d_z$. Noting that ${\rm d}_z \in {\cal T}_{\K}(z)$ is characterized by $\nabla (\psi \circ \lambda)(A)^*d_A -d_s\leq 0$, we obtain that (\ref{meq:projD-2}) is satisfied for this case.
For any $z'=(A',s')$, denote
$$
(X(z'),t(z'))=\Pi_{{\cal K}}(z') \mbox{ and } \mu (z')=t(z')-s',
$$
then $t(z')= (\psi\circ \lambda)(X(z'))=N(X(z'))$ when $z' \notin {\rm int}\,{\cal K}$.
For the case when $z=(A,s) \in {\rm bdry} \K \setminus\{0\}$, let  $d_z \notin {\cal T}_{\K}(z)$, then $\nabla (\psi \circ \lambda)(A)^*d_A -d_s\geq 0$ and
$$
z+ \gamma d_z \in   {\rm int }[\mathbb S^m \times \mathbb R \setminus (\K\cup[\K]^{\circ})]
$$
for small $ \gamma>0$. From the definition of $(X(\cdot),\mu(\cdot))$, we have
$$
\begin{array}{l}
\mu(z+\gamma d_z) \nabla (\psi\circ \lambda)(X(z+\gamma d_z))+X(z+\gamma d_z)-(A+\gamma d_A)=0\\[3pt]
(\psi\circ \lambda)(X(z+\gamma d_z))-\mu(z+\gamma d_z) -(s+\gamma d_s)=0.
\end{array}
$$
Thus we obtain
$$
[(\psi\circ\lambda)(X(z+\gamma d_z)) -(s+\gamma d_s)] \nabla (\psi\circ \lambda) (X(z+\gamma d_z))+X(z+\gamma d_z)-(A+\gamma d_A)=0.
$$
Noting $A=X(z), s=(\psi\circ\lambda)(A)$, we obtain from the above relation that
\begin{equation}\label{meq:Ah1}
\begin{array}{r}
[(\psi\circ \lambda)(X(z+\gamma d_z)) -(\psi\circ \lambda)(X(z))-\gamma d_s] \nabla (\psi \circ \lambda)(X(z+\gamma d_z))\\[6pt]
\quad \quad +X(z+\gamma d_z)-(A+\gamma d_A)=0.
\end{array}
\end{equation}
From the mean-value theorem, there exists $\beta_{\gamma} \in (0,1)$ such that
\begin{equation}\label{meq:Ah2}
\begin{array}{l}
(\psi\circ \lambda)(X(z+\gamma d_z)) -(\psi\circ \lambda)(X(z))\\[8pt]
=\nabla (\psi\circ \lambda)( (1-\beta_{\gamma})X(z)+
\beta_{\gamma}X(z+\gamma d_z))^*(X(z+\gamma d_z)-X(z)).
\end{array}
\end{equation}
Obviously we have $\beta_{\gamma} \rightarrow 0$ when $\gamma \searrow 0$. Combining (\ref{meq:Ah1}) and (\ref{meq:Ah2}), we obtain
$$
\begin{array}{l}
\left[{\cal I}+ \nabla (\psi\circ \lambda) (X(z+\gamma d_z)) \nabla (\psi\circ \lambda)( (1-\beta_{\gamma})X(z)+
\beta_{\gamma}X(z+\gamma d_z))^*\right]\displaystyle \frac{X(z+\gamma d_z)-X(z)}{\gamma}\\[6pt]
=\nabla(\psi\circ \lambda) (X(z+\gamma d_z)) d_s+d_A.
\end{array}
$$
When $\gamma >0$ is small enough, then the operator
$$\left[{\cal I}+ \nabla (\psi\circ \lambda) (X(z+\gamma d_z)) \nabla (\psi\circ \lambda)( (1-\beta_{\gamma})X(z)+
\beta_{\gamma}X(z+\gamma d_z))^*\right]
$$
is nonsingular, and thus
$$
\begin{array}{l}
\displaystyle \frac{X(z+\gamma d_z)-X(z)}{\gamma}=\\[6pt]
\left[{\cal I}+ \nabla (\psi\circ \lambda) (X(z+\gamma d_z)) \nabla (\psi\circ \lambda)( (1-\beta_{\gamma})X(z)+
\beta_{\gamma}X(z+\gamma d_z))^*\right]^{-1}(\nabla(\psi\circ \lambda) (X(z+\gamma d_z)) d_s+d_A).
\end{array}
$$
Since the right-hand side of the above equation has the limit when $\gamma \searrow 0$, we have that $X(\cdot)$ is directionally differentiable at $z$ along $d_z$, and
$$
\begin{array}{ll}
X'(z;d_z)&=[{\cal I}+ \nabla (\psi\circ \lambda)(A)\nabla (\psi\circ \lambda)(A)^*]^{-1}(\nabla (\psi\circ \lambda)(A) d_s+d_A)\\[6pt]
& =\left[{\cal I}-\displaystyle \frac{\nabla (\psi\circ \lambda)(A)\nabla (\psi\circ \lambda)(A)^*}{1+\|(\psi\circ \lambda)(A)\|_F^2} \right](\nabla (\psi\circ \lambda)(A) d_s+d_A)\\[12pt]
&=d_A+\displaystyle \frac{d_s-\nabla (\psi\circ \lambda)(A)^*d_A}{1+\|\nabla (\psi\circ \lambda)(A)\|^2_F}\nabla (\psi\circ \lambda)(A).
\end{array}
$$
Noting $t(z)=(\psi\circ\lambda )(X(z))$, we obtain from the above expression of $X'(z,d_z)$ that
$$
\begin{array}{ll}
t'(z;d_z) &= \nabla (\psi\circ \lambda)(A)^* X'(z;d_z)\\[6pt]
&=d_s-\displaystyle \frac{d_s-\nabla (\psi\circ \lambda)(A)^*d_A}{1+\|\nabla (\psi\circ \lambda)(A)\|^2_F}.
\end{array}
$$
This implies the expression of (\ref{meq:projD-2}).

Now we prove (v).  From (iv), as $z=(A,s) \in {\rm bdry} [\K]^{\circ} \setminus\{0\}$, we obtain
\begin{equation}\label{eq:qDir}
\begin{array}{l}
\Pi'_{[\K]^{\circ}}(z;d_z)=
\left[
\begin{array}{c}
d_A\\[14pt]
d_s
\end{array}
\right]-\displaystyle \frac{[\nabla (\psi_*\circ \lambda)(A)^*d_A +d_s]_+}{1+\|\nabla(\psi_*\circ \lambda)(A)\|^2_F}\left[
\begin{array}{c}
\nabla (\psi_*\circ \lambda)(A)\\[14pt]
1
\end{array}
\right].
\end{array}
\end{equation}
Noting the following identity
$$
\Pi_{\K}(z)=z-\Pi_{[\K]^{\circ}}(z)
$$
and $\Pi'_{\K}(z;d_z)=d_z-\Pi'_{[\K]^{\circ}}(z;d_z)$, we obtain (v) from (\ref{eq:qDir}).
The proof is completed. \hfill $\Box$\\
Noting that $N=N_*$ when $\psi$ is the $l_2$-norm of $\mathbb R^m$,  $N$ is twice   differentiable  norms,  we have the following results about the directional derivatives of the projection operator onto the second-order cone.
\begin{corollary}\label{mCorPiDirction2}
For the projection operator onto the second-order cone ${\cal K}_2$, we have
\begin{itemize}
\item[{\rm (i)}] If $z=(A,s) \in {\rm int } \K_2$, then $\Pi'_{\K_2}(z;d_z)=d_z$;
\item[{\rm (ii)}] If $z=(A,s) \in {\rm int } [\K_2]^{\circ}$, then $\Pi'_{\K_2}(z;d_z)=0$;
\item[{\rm (iii)}] If $z=(A,s) \in {\rm int }[\mathbb S^m \times \mathbb R\setminus (\K_2\cup[\K_2]^{\circ})]$, then
\begin{equation}\label{meq:projD-12Cor}
\begin{array}{l}
\Pi'_{\K_2}(z;d_z)=\displaystyle \frac{1}{2}\left[
\begin{array}{rc}
{\cal I}+\displaystyle \frac{s}{\|A\|_F}{\cal I}-\displaystyle \frac{s}{\|A\|_F}\displaystyle \frac{AA^*}{\|A\|_F^2} & \displaystyle \frac{A}{\|A\|_F}\\[14pt]
\displaystyle \frac{A^*}{\|A\|_F} & 1
\end{array}
\right]\left
[
\begin{array}{l}
d_A\\[18pt]
d_s
\end{array}
\right].
\end{array}
\end{equation}
\item[{\rm (iv)}] If $z=(A,s) \in {\rm bdry} \K_2 \setminus\{0\}$, then
  \begin{equation}\label{meq:projD-2Cor}
\begin{array}{l}
\Pi'_{\K_2}(z;d_z)=
\left[
\begin{array}{c}
d_A\\[14pt]
d_s
\end{array}
\right]-\displaystyle \frac{1}{2}\left[\displaystyle\frac{A^*d_A}{\|A\|_F} -d_s\right]_+\cdot\left[
\begin{array}{c}
\displaystyle \frac{A}{\|A\|_F}\\[14pt]
-1
\end{array}
\right].
\end{array}
\end{equation}
\item[{\rm (v)}] If $z=(A,s) \in {\rm bdry} [\K_2]^{\circ} \setminus\{0\}$, then
  \begin{equation}\label{eq:projD-3}
\begin{array}{l}
\Pi'_{\K_2}(z;d_z)=
\displaystyle \frac{1}{2}\left[\displaystyle\frac{A^*d_A}{\|A\|_F} +d_s\right]_+\cdot\left[
\begin{array}{c}
\displaystyle \frac{A}{\|A\|_F}\\[14pt]
1
\end{array}
\right].
\end{array}
\end{equation}
\item[{\rm (vi)}] If $z=(A,s)=(0,0)$, then $\Pi'_{\K_2}(z;d_z)=\Pi_{\K_2}({\rm d}_z)$.
\end{itemize}
\end{corollary}
{\bf Proof}. Conclusions (i),(ii),(iii) and (vi) come from Corollary \ref{coromPiD2}. Noting
$$
\nabla (\psi\circ \lambda)(A)=\displaystyle \frac{A}{\|A\|_F},\,\, \|\nabla (\psi\circ \lambda)(A)\|_F=1,
$$
we obtain (iv) and (v) from Theorem \ref{mthPiDirction} (iv)(v). \hfill $\Box$\\
Let
$$
\bar G=\left[
\begin{array}{ll}
G(z) & 0\\[14pt]
0 & 1
\end{array}
\right],\, v(z)=\left[
\begin{array}{c}
\nabla (\psi\circ \lambda)(X(z))\\[14pt]
-1
\end{array}
\right].
$$
Then when $z=(A,s) \in {\rm int }[\mathbb R^{m+1}\setminus (\K\cup[\K]^{\circ})]$, we have from (\ref{meq:projD}) that
\begin{equation}\label{meq:ABS}
{\rm D}\Pi_{\K}(z)=\bar G(z)^{-1}-\displaystyle \frac{\bar G(z)^{-1}v(z)v(z)^*\bar G(z)^{-1}}{v(z)^*\bar G_N(z)^{-1}v(z)}.
\end{equation}
Let ${\cal I}': {\cal H}\rightarrow {\cal H}$ be the identity mapping, we usually denote it by
$${\cal I}'=\left[
    \begin{array}{ll}
    {\cal I} &0\\[4pt]
    0 & 1
    \end{array}
    \right],
    $$
where ${\cal I}: \mathbb{S}^{m} \rightarrow \mathbb{S}^{m} $ is the identity mapping.

B-subdifferential of the projection operator onto $\K$ is also very important in stability of Problem (\ref{Pconic}), especially the strong regularity of Kurash-Kuhn-Tucker system of the problem.
\begin{theorem}\label{mthBsubdiff}
Let $\psi: \mathbb R^m \rightarrow \mathbb R_+$ be a  differentiable symmetric gauge function. Let $z=(A,s)\in {\cal H}$ be a given point with $A=P{\rm Diag}(\lambda (A))P^T$, then
  the following results hold:
\begin{itemize}
\item[{\rm (i)}] If $z=(A,s) \in {\rm int } \K$, then $\partial_B \Pi_{\K}(z)=\{{\cal I}'\}$,
\item[{\rm (ii)}] If $z=(A,s) \in {\rm int } [\K]^{\circ}$, then $\partial_B \Pi_{\K}(z)=\{0\}$;
\item[{\rm (iii)}] If $z=(A,s) \in {\rm int }[{\cal H}\setminus (\K\cup[\K]^{\circ})]$
and
$\psi$ is a twice symmetric gauge function,  then
$\partial_B \Pi_{\K}(z)=\{{\rm D}\Pi_{\K}(z)\}$, where ${\rm D}\Pi_{\K}(z)$ is calculated by
(\ref{eq:mprojD});
\item[{\rm (iv)}] If $z=(A,s) \in {\rm bdry} \K \setminus\{0\}$ and
$\psi$ is a twice symmetric gauge function, then
  \begin{equation}\label{meq:projD-2Bd}
  \partial_B \Pi_{\K}(z)=\left\{{\cal I}', \left[\begin{array}{ll}
  {\cal I} & 0\\[6pt]
  0 & 1
  \end{array}
  \right]-\left
  [
  \begin{array}{c}
 \nabla (\psi\circ \lambda)(A)\\[6pt]
  -1
  \end{array}
  \right]
  \displaystyle \frac{[\nabla (\psi\circ \lambda)(A)^* \,\, -1]}{ 1+\|\nabla (\psi\circ \lambda)(A)\|^2_F} \right\}.
  \end{equation}
\item[{\rm (v)}] If $z=(A,s) \in {\rm bdry} [\K]^{\circ} \setminus\{0\}$
 and
$\psi_*$ is a twice symmetric gauge function, then
  \begin{equation}\label{meq:projD-3Bd}
\begin{array}{l}
\partial_B \Pi_{\K}(z)=\left\{0,
\left[
\begin{array}{c}
\nabla (\psi_*\circ \lambda)(A)\\[14pt]
1
\end{array}
\right]\displaystyle \frac{[\nabla (\psi_*\circ \lambda)(A)^*\,\, 1]}{1+\|\nabla (\psi_*\circ \lambda)(A)\|^2_F}\right\}.
\end{array}
\end{equation}
\item[{\rm (vi)}] If $z=(A,s)=(0,0)$ and
$\psi$ is a twice symmetric gauge function, then\footnote{In the following formula, $\limsup_{z \stackrel{{\cal D}}\rightarrow 0}A(z)$ for some mapping $A:\mathbb S^{m}\times \mathbb R\rightarrow \mathbb S^{m}\times \mathbb R$ is defined by
    $$
    \limsup_{z \stackrel{{\cal D}}\rightarrow 0}A(z)=\left\{A \in \mathbb S^{m}\times \mathbb R: \exists z^k \in {\cal D}, z^k \rightarrow 0\mbox{ such that } A=\lim_{k \rightarrow \infty} A(z^k)\right\}.
    $$
    }
\begin{equation}\label{eq:0Bd}
\partial_B \Pi_{\K}(z)=\left\{0,{\cal I}'\right\}\cup \limsup_{z \stackrel{{\cal D}}\rightarrow 0}\left\{\bar G(z)^{-1}-\displaystyle \frac{\bar G(z)^{-1}v(z)v(z)^*\bar G(z)^{-1}}{v(z)^*\bar G(z)^{-1}v(z)}\right\},
\end{equation}
where ${\cal D}={\rm int }[\mathbb S^{m}\times \mathbb R\setminus (\K\cup[\K]^{\circ})]$.
\end{itemize}
\end{theorem}
{\bf Proof}. In cases (i),(ii) and (iii), $\Pi_{\K}$ is differentiable at $z$, the results are from Proposition \ref{mpropPiD}.  For case (iv), $z=(A,s) \in {\rm bdry} \K \setminus\{0\}$, for  constructing the B-subdifferential of
 $\Pi_{\K}$ at $z$,
we have two ways of $z'$ approaches $z$, one is $z' \stackrel{{\rm int } \K} \longrightarrow z$ and the other is $z' \stackrel{{\cal D}}\rightarrow z$, where ${\cal D}={\rm int }[\mathbb S^{m}\times \mathbb R\setminus (\K\cup[\K]^{\circ})]$.
For $z' \stackrel{{\cal D}}\rightarrow z$, we obtain ${\cal I}'\in \partial_B \Pi_{\K}(z)$. For $z' \stackrel{{\cal D}}\rightarrow z$, one has
$$
\lim_{z' \stackrel{{\cal D}}\rightarrow z} X(z')=A, \, \lim_{z' \stackrel{{\cal D}}\rightarrow z} \mu(z')=0,
$$
which implying
$$
\lim_{z' \stackrel{{\cal D}}\rightarrow z} G(z')={\cal I},\, \lim_{z' \stackrel{{\cal D}}\rightarrow z}G(z)^{-1}\nabla (\psi \circ \lambda)(X(z))=\nabla (\psi\circ \lambda)(A).
$$
Thus we obtain
$$
\lim_{z' \stackrel{{\cal D}}\rightarrow z} \beta(z')=(1+\|\nabla (\psi\circ \lambda)(A)\|^2_F)^{-1}
$$
and
$$
\begin{array}{l}
\lim_{z' \stackrel{{\cal D}}\rightarrow z}{\rm D}\Pi_{\K}(z')\\[14pt]
=\left[
\begin{array}{ll}
{\cal I} & 0\\[14pt]
0 & 1
\end{array}
\right] -(1+\|\nabla (\psi \circ \lambda) (A)\|^2_F)^{-1}

\left[
\begin{array}{c}
\nabla (\psi \circ \lambda) (A)\\[14pt]
-1
\end{array}
\right][\nabla (\psi \circ \lambda) (A)^* \,\,\, -1].
\end{array}
$$
Combining the above two limits, we obtain (\ref{meq:projD-2Bd}).

For case (v), $z=(A,s) \in {\rm bdry} [\K]^{\circ} \setminus\{0\}$, for  constructing the B-subdifferential of
 $\Pi_{\K}$ at $z$,
we have two ways of $z'$ approaches $z$, one is $z' \stackrel{[\K]^{\circ}} \longrightarrow z$ and the other is $z' \stackrel{{\cal D}}\rightarrow z$, where ${\cal D}={\rm int }[\mathbb S^{m}\times \mathbb R\setminus (\K\cup[\K]^{\circ})]$.
For $z' \stackrel{{\cal D}}\rightarrow z$, we obtain $0\in \partial_B \Pi_{\K}(z)$. For $z' \in {\cal D}$, one has
$$
{\rm D}\Pi_{\K}(z')={\cal I}'-{\rm D}\Pi_{[\K]^{\circ}}(z').
$$
 From Proposition \ref{mpropqPiDq}, for $z' \in {\cal D}$,

 \begin{equation}\label{meq:projDqprime}
\begin{array}{l}
{\rm D}\Pi_{\K^{\circ}}(z') =
\left[
\begin{array}{ll}
G_*(z')^{-1} & 0\\[8pt]
0 & 1
\end{array}
\right]  \\[4pt]
\quad -\beta_*(z')\left[
\begin{array}{c}
G_*(z')^{-1}\nabla (\psi_*\circ \lambda)(X(z'))\\[8pt]
1
\end{array}
\right][\nabla (\psi_*\circ \lambda)(X(z'))^*G_*(z')^{-1} \,\,\, 1],
\end{array}
\end{equation}
where
$$\beta_*(z')=(1+\nabla (\psi_*\circ \lambda)(X(z'))^*G_*(z')^{-1}\nabla (\psi_*\circ \lambda)(X(z')))^{-1},\, G_*(z')={\cal I}+\mu(z') \nabla^2 (\psi_*\circ \lambda)(X(z')),$$
 where $X(z')=P[{\rm Diag}u_*(z')]P^T$ with $(u_*(z'),\mu_*(z'))\in \mathbb R^m \times \mathbb R_{++}$ being a solution of the system of equations $F_*(u,\mu;z')=0$. For $z' \stackrel{{\cal D}}\rightarrow z$, one has
$$
\lim_{z' \stackrel{{\cal D}}\rightarrow z} u_*(z')=\lambda (A), \, \lim_{z' \stackrel{{\cal D}}\rightarrow z} \mu_*(z')=0,
$$
which implying
$$
\lim_{z' \stackrel{{\cal D}}\rightarrow z} G_*(z')={\cal I},\, \lim_{z' \stackrel{{\cal D}}\rightarrow z}G_*(z')^{-1}\nabla (\psi_* \circ \lambda)(X(z'))=\nabla (\psi_* \circ \lambda) (A).
$$
 Thus we obtain
$$
\lim_{z' \stackrel{{\cal D}}\rightarrow z} \beta_{*}(z')=(1+\|\nabla (\psi_* \circ \lambda)(A)\|^2_F)^{-1}
$$
and
$$
\begin{array}{l}
\displaystyle \lim_{z' \stackrel{{\cal D}}\rightarrow z}{\rm D}\Pi_{\K}(z')={\cal I}'-\lim_{z' \stackrel{{\cal D}}\rightarrow z}{\rm D}\Pi_{[\K]^{\circ}}(z')\\[14pt]
=(1+\|\nabla (\psi_* \circ \lambda)(A)\|^2_F)^{-1}
\left[
\begin{array}{c}
\nabla (\psi_* \circ \lambda)(A)\\[14pt]
1
\end{array}
\right][\nabla (\psi_* \circ \lambda)(A)^* \,\,\, 1].
\end{array}
$$
Combining the above two limits, we obtain (\ref{meq:projD-3Bd}).

When $z=(A,s)=(0,0)$,  for  constructing the B-subdifferential of
 $\Pi_{\K}$ at $0$,
we have three  ways of $z'$ approaches $z$:
 $$
   z' \stackrel{\K} \longrightarrow 0,\, z' \stackrel{[\K]^{\circ}} \longrightarrow 0 \mbox{ and } z' \stackrel{{\cal D}}\rightarrow 0,$$
   where ${\cal D}={\rm int }[\mathbb S^{m}\times \mathbb R\setminus (\K\cup[\K]^{\circ})]$. The same analysis as in (iv) or (v) yields the formula $\partial_B \Pi_{\K}(0)$ in (\ref{eq:0Bd}). The proof is completed. \hfill $\Box$\\
 For the second-order cone, we have the following results about the B-subdifferential of $\Pi_{\K_2}$.
\begin{corollary}\label{Cor2Bsubdiff}
For the  second-order cone ${\cal K}_2$ in $\mathbb S^{m}\times \mathbb R$, we have
\begin{itemize}
\item[{\rm (i)}] If $z=(A,s) \in {\rm int } \K_2$, then $\partial_B \Pi_{\K_2}(z)=\{{\cal I}'\}$;
\item[{\rm (ii)}] If $z=(A,s) \in {\rm int } [\K_2]^{\circ}$, then $\partial_B \Pi_{\K_2}(z)=\{0\}$;
\item[{\rm (iii)}] If $z=(A,s) \in {\rm int }[\mathbb S^{m}\times \mathbb R\setminus (\K_2\cup[\K_2]^{\circ})]$, then
$$\partial_B \Pi_{\K_2}(z)=\left\{\displaystyle \frac{1}{2}\left[
\begin{array}{rc}
{\cal I}+\displaystyle \frac{s}{\|A\|_F}{\cal I}-\displaystyle \frac{s}{\|A\|_F}\displaystyle \frac{AA^*}{\|A\|_F^2} & \displaystyle \frac{A}{\|A\|_F}\\[14pt]
\displaystyle \frac{A^*}{\|A\|_F} & 1
\end{array}
\right]\right\};$$
\item[{\rm (iv)}] If $z=(A,s) \in {\rm bdry} \K_2 \setminus\{0\}$, then
    \begin{equation}\label{meq:corojD-2Bd}
  \partial_B \Pi_{\K_2}(z)=\left\{{\cal I}', \displaystyle \frac{1}{2}\left[\begin{array}{cc}
  2{\cal I} -\displaystyle \frac{AA^*}{\|A\|_F^2} & \displaystyle \frac{A}{\|A\|_F}\\[12pt]
  \displaystyle \frac{A^*}{\|A\|_F} & 1
  \end{array}
  \right] \right\}.
  \end{equation}
\item[{\rm (v)}] If $z=(A,s) \in {\rm bdry} [\K_2]^{\circ} \setminus\{0\}$, then
  \begin{equation}\label{meq:coroojD-3Bd}
\begin{array}{l}
\partial_B \Pi_{\K_2}(z)=\left\{0,\,
 \displaystyle \frac{1}{2}\left[\begin{array}{cc}
  \displaystyle \frac{AA^*}{\|A\|_F^2} & \displaystyle \frac{A}{\|A\|_F}\\[12pt]
  \displaystyle \frac{A^*}{\|A\|_F} & 1
  \end{array}
  \right]\right\}.
\end{array}
\end{equation}
\item[{\rm (vi)}] If $z=(A,s)=(0,0)$,
\begin{equation}\label{meq:coro0Bd}
\partial_B \Pi_{\K_2}(0)=\left\{ 0,{\cal I}'\right\}\cup \left\{ \displaystyle \frac{1}{2}\left[
\begin{array}{cc}
2 a({\cal I}-YY^*)+YY^* & Y\\[12pt]
Y^*  & 1
\end{array}
\right]:a \in [0,1],\|Y\|_F=1
\right\}.
\end{equation}
\end{itemize}
\end{corollary}
{\bf Proof}. Assertions (i), (ii) and (iii) come from Corollary \ref{coromPiD2}. For (iv),  $z=(A,s) \in {\rm bdry} \K_2 \setminus\{0\}$, $A \ne 0$, $\nabla (\phi \circ \lambda)(A)=A/\|A\|_F$ and $\|\nabla (\phi \circ \lambda)(A)\|_F=1$.  The formula (\ref{meq:corojD-2Bd}) comes from (\ref{meq:projD-2Bd}). For $\psi$ being the $l_2$-norm, $\psi=\psi_*$ and in this case $\|\nabla (\phi_* \circ \lambda)(A)\|_2=1$ and (\ref{meq:coroojD-3Bd}) comes from (\ref{meq:projD-3Bd}).

 For ${\cal D}={\rm int }[\mathbb S^{m}\times \mathbb R\setminus (\K_2\cup[\K_2]^{\circ})]$, if $z'=(w',s') \in {\cal D}$, then we have from Corollary \ref{coromPiD2} (iii) that
\begin{equation}\label{meq:projD2c}
\begin{array}{ll}
{\rm D}\Pi_{\K_2}(z')& =\displaystyle \frac{1}{2}\left[
\begin{array}{rc}
{\cal I}+\displaystyle \frac{s'}{\|A'\|_F}{\cal I}-\displaystyle \frac{s'}{\|A'\|_F}\displaystyle \frac{A'{A'}^*}{\|A'\|_F^2} & \displaystyle \frac{A'}{\|A'\|_F}\\[14pt]
\displaystyle \frac{{A'}^*}{\|A'\|_F} & 1
\end{array}
\right]\\[18pt]
&=\displaystyle \frac{1}{2}\left[
\begin{array}{rc}
\left(1+\displaystyle \frac{s'}{\|A'\|_F}\right)\left({\cal I}-\displaystyle \frac{A'{A'}^*}{\|A'\|_F^2} \right)& \displaystyle \frac{A'}{\|A'\|_F}\\[14pt]
\displaystyle \frac{{A'}^*}{\|A'\|_F} & 1
\end{array}
\right].
\end{array}
\end{equation}
Noting that when $z' \stackrel{{\cal D}}\rightarrow z$, the outer limit of $$\left\{\displaystyle \frac{s'}{\|A'\|_F},\displaystyle \frac{A'}{\|A'\|_F}\right\}$$
coincides with $\{(a, Y): a\in [0,1], \|Y\|_F=1\}$, and we obtain the formula (\ref{meq:coro0Bd}). The proof is completed. \hfill $\Box$
\section{Calculating $\nabla^2 (\psi\circ \lambda)$ and $G(z)^{-1}$}\label{Section-computing}
\setcounter{equation}{0}
In order to simplify formulas for directional derivatives and B-subdifferentials of $\Pi_{{\cal K}}$ when $\psi$ is twice continuously differentiable, we develop formulas for $\nabla^2 (\psi\circ \lambda)(A)$ and $G(z)^{-1}H$.
\begin{proposition}\label{prop:DF}
Let  $A$ have the spectral decomposition $A=P{\rm Diag}\lambda (A) P^T$, $w_1,\ldots, w_r$ be $r$ distinct values of $m$ eigenvalues of $A$ with
$$
w_1=\lambda_1(A)=\cdots=\lambda_{k_1}(A),w_2=\lambda_{k_1+1}(A)=\cdots=\lambda_{k_2}(A),
w_r=\lambda_{k_{r-1}+1}(A)=\cdots=\lambda_{m}(A).$$
  Denote
$$
\alpha_1=\{1,\ldots, k_1\}, \alpha_2=\{k_1+1,\ldots, k_2\}, \ldots, \alpha_r=\{k_{r-1}+1,\ldots, k_r\}.
$$
Let $\psi: \mathbb R^m \rightarrow \mathbb R$ be a symmetric function, twice differentiable at the point $w \in \mathbb R^m$,
and $P$ be a permutation matrix such that $Pw=w$. Then
\begin{itemize}
\item[(1)]there are real numbers $b_1,\ldots, b_r$ and a symmetric matrix $(a_{ij})_{i,j=1}^r$ such that
\begin{equation}\label{2ndpsimu}
\begin{array}{ll}
\nabla^2 \psi (w)=&\left[
\begin{array}{cccc}
b_1I_{|\alpha_1|} & & &\\[4pt]
& b_2I_{|\alpha_2|}& &\\[4pt]
& & \ddots  &\\[4pt]
& & & b_rI_{|\alpha_r|}
\end{array}
\right]\\[16pt]
& +\left [
\begin{array}{cccc}
a_{11}\textbf{1}_{|\alpha_1|}\textbf{1}_{|\alpha_1|}^T &a_{12}\textbf{1}_{|\alpha_1|}\textbf{1}_{|\alpha_2|}^T&\cdots & a_{1r}\textbf{1}_{|\alpha_1|}\textbf{1}_{|\alpha_r|}^T\\[10pt]
a_{21}\textbf{1}_{|\alpha_2|}\textbf{1}_{|\alpha_1|}^T&a_{22}\textbf{1}_{|\alpha_2|}\textbf{1}_{|\alpha_2|}^T &\cdots & a_{2r}\textbf{1}_{|\alpha_2|}\textbf{1}_{|\alpha_r|}^T\\[10pt]
\vdots & \vdots &\vdots& \vdots\\[10pt]
a_{r1}\textbf{1}_{|\alpha_r|}\textbf{1}_{|\alpha_1|}^T&a_{r2}\textbf{1}_{|\alpha_r|}\textbf{1}_{|\alpha_2|}^T & \cdots & a_{rr}\textbf{1}_{|\alpha_r|}\textbf{1}_{|\alpha_r|}^T
\end{array}
\right].
\end{array}
\end{equation}
\item[(2)]there are real numbers $c_1,\ldots, c_r$ such that
$$
\nabla \psi (w)=\left[
\begin{array}{c}
c_1\textbf{1}_{|\alpha_1|}\\[2pt]
\vdots\\[2pt]
c_r\textbf{1}_{|\alpha_r|}
\end{array}
\right],
$$
and for any $H \in \mathbb S^m$,
\begin{equation}\label{eq:2ndH}
\begin{array}{l}
\nabla^2 (\psi\circ\lambda)(A)[H]=\\[10pt]
P\left[
\begin{array}{ccc}
\left(\displaystyle \sum_{j=1}^ra_{1j} {\rm Tr}\left(P_{\alpha_j}^THP_{\alpha_j}\right)\right)I_{|\alpha_1|}  & &\\[4pt]
&  \ddots  &\\[4pt]
 & & \left(\displaystyle \sum_{j=1}^ra_{rj} {\rm Tr}\left(P_{\alpha_j}^THP_{\alpha_j}\right)\right) I_{|\alpha_r|}
\end{array}
\right]P^T\\[12pt]
+P\left[
\begin{array}{ccc}
\delta_{11}P_{\alpha_1}^THP_{\alpha_1}  & \cdots & \delta_{1r}P_{\alpha_1}^THP_{\alpha_r}\\[4pt]
\vdots &  \ddots  & \vdots\\[4pt]
\delta_{r1}P_{\alpha_r}^THP_{\alpha_1}  & \cdots & \delta_{rr}P_{\alpha_r}^THP_{\alpha_r}
\end{array}
\right]P^T
\end{array}
\end{equation}
and
\begin{equation}\label{eq:2ndDSF}
\begin{array}{l}
\nabla^2 (\psi\circ\lambda)(A)[H,H]\\[6pt]
 =\displaystyle \sum_{i=1}^r \sum_{j=1}^r a_{ij}{\rm Tr}
(P_{\alpha_i}^THP_{\alpha_i}){\rm Tr}
(P_{\alpha_j}^THP_{\alpha_j})
+\displaystyle \sum_{i=1}^r \sum_{j=1}^r \delta_{ij}\displaystyle \sum_{i'\in \alpha_i,j'\in \alpha_j}(p_{i'}^THp_{j'})^2,
\end{array}
\end{equation}
\end{itemize}
where
$$
\left\{
\begin{array}{l}
\delta_{ii}=b_i, i=1,\ldots, r,\\[8pt]
\delta_{ij}=\displaystyle \frac{c_i-c_j}{w_i-w_j},i,j=1,\ldots r, i \ne j.
\end{array}
\right.
$$
\end{proposition}
{\bf Proof}. The formula (\ref{2ndpsimu}) is an equivalent version of the formula in Lemma \ref{lem:LS2.1}(ii).
In view of Lemma \ref{lem:LS2.1}, under the assumptions here we can express ${\cal A}$ as
\begin{equation}\label{2ndpsimua}
\begin{array}{l}
{\cal A}=-\left[
\begin{array}{ccc}
\delta_{11}I_{|\alpha_1|}  & &\\[4pt]
 & \ddots  &\\[4pt]
& & \delta_{rr} I_{|\alpha_r|}
\end{array}
\right] +\left [
\begin{array}{ccc}
\delta_{11}\textbf{1}_{|\alpha_1|}\textbf{1}_{|\alpha_1|}^T  &\cdots & \delta_{1r}\textbf{1}_{|\alpha_1|}\textbf{1}_{|\alpha_r|}^T\\[10pt]
 \vdots &\vdots& \vdots\\[10pt]
\delta_{r1}\textbf{1}_{|\alpha_r|}\textbf{1}_{|\alpha_1|}^T & \cdots & \delta_{rr}\textbf{1}_{|\alpha_r|}\textbf{1}_{|\alpha_r|}^T
\end{array}
\right].
\end{array}
\end{equation}
For $\widehat H=P^THP$, define
$$
h_1=\left(
\begin{array}{c}
p_1^THp_1\\[2pt]
\vdots\\[2pt]
p_{|\alpha_1|}^THp_{|\alpha_1|}
\end{array}
\right), \, h_2=\left(
\begin{array}{c}
p_{|\alpha_1|+1}^THp_{|\alpha_1|+1}\\[2pt]
\vdots\\[2pt]
p_{|\alpha_2|}^THp_{|\alpha_2|}
\end{array}
\right),\ldots, h_r=\left(
\begin{array}{c}
p_{|\alpha_{r-1}|+1}^THp_{|\alpha_{r-1}|+1}\\[2pt]
\vdots\\[2pt]
p_{|\alpha_r|}^THp_{|\alpha_r|}
\end{array}
\right),
$$
we have from (\ref{2ndpsimu}) and (\ref{2ndpsimua}) that
$$
\begin{array}{l}
\nabla^2 (\psi\circ\lambda)(A)[H]=P\left({\rm Diag}(\nabla^2 \psi(\lambda(A)){\rm diag}\,\widehat H)+{\cal A}\circ \widehat H\right)P^T\\[10pt]
=P\left( {\rm Diag}\left(\left[ \begin{array}{c}
b_1h_1\\[2pt]
\vdots\\[2pt]
b_rh_r
\end{array} \right]+
\left[
\begin{array}{c}
\left(a_{11}\textbf{1}_{|\alpha_1|}^Th_1+a_{12}\textbf{1}_{|\alpha_2|}^Th_2+\cdots+a_{1r}\textbf{1}_{|\alpha_r|}^Th_r \right)\textbf{1}_{|\alpha_1|}\\[4pt]
\vdots\\[4pt]
\left(a_{r1}\textbf{1}_{|\alpha_1|}^Th_1+a_{r2}\textbf{1}_{|\alpha_2|}^Th_2+\cdots+a_{rr}\textbf{1}_{|\alpha_r|}^Th_r \right)\textbf{1}_{|\alpha_r|}
\end{array}
\right]\right)
\right)P^T\\[10pt]
\quad -P\left[
\begin{array}{cccc}
\delta_{11}I_{|\alpha_1|}\circ P_{\alpha_1}^THP_{\alpha_1} & & &\\[8pt]
& \delta_{22}I_{|\alpha_2|}\circ P_{\alpha_2}^THP_{\alpha_2}  & &\\[8pt]
& &\ddots &\\[8pt]
&&&\delta_{rr}I_{|\alpha_r|}\circ P_{\alpha_r}^THP_{\alpha_r}
\end{array}
\right]P^T\\[10pt]
\quad +P\left[
\begin{array}{cccc}
\delta_{11}P_{\alpha_1}^THP_{\alpha_1} &\delta_{12}P_{\alpha_1}^THP_{\alpha_2} & \cdots & \delta_{1r}P_{\alpha_1}^THP_{\alpha_r}\\[8pt]
\delta_{21}P_{\alpha_2}^THP_{\alpha_1} &\delta_{22}P_{\alpha_2}^THP_{\alpha_2} & \cdots & \delta_{2r}P_{\alpha_2}^THP_{\alpha_r}\\[8pt]
\vdots &\vdots &\vdots & \vdots\\[8pt]
\delta_{r1}P_{\alpha_r}^THP_{\alpha_1} &\delta_{r2}P_{\alpha_r}^THP_{\alpha_2} & \cdots & \delta_{rr}P_{\alpha_r}^THP_{\alpha_r}
\end{array}
\right]P^T\\[10pt]
=P\left[
\begin{array}{ccc}
\left(\displaystyle \sum_{j=1}^ra_{1j} {\rm Tr}\left(P_{\alpha_j}^THP_{\alpha_j}\right)\right)I_{|\alpha_1|}  & &\\[4pt]
&  \ddots  &\\[4pt]
 & & \left(\displaystyle \sum_{j=1}^ra_{rj} {\rm Tr}\left(P_{\alpha_j}^THP_{\alpha_j}\right)\right) I_{|\alpha_r|}
\end{array}
\right]P^T\\[12pt]
\quad +P\left[
\begin{array}{ccc}
\delta_{11}P_{\alpha_1}^THP_{\alpha_1}  & \cdots & \delta_{1r}P_{\alpha_1}^THP_{\alpha_r}\\[4pt]
\vdots &  \ddots  & \vdots\\[4pt]
\delta_{r1}P_{\alpha_r}^THP_{\alpha_1}  & \cdots & \delta_{rr}P_{\alpha_r}^THP_{\alpha_r}
\end{array}
\right]P^T,
\end{array}
$$
which verifies (\ref{eq:2ndH}).

The second-order directional derivative of $(\psi\circ\lambda)$ at $A$ along $H$ is
$$
\begin{array}{ll}
\nabla^2 (\psi\circ\lambda)(A)[H,H]
& =\nabla^2 \psi (w)[{\rm diag}\widehat H,{\rm diag}\widehat H]+\langle {\cal A}, \widehat H \circ \widehat H\rangle\\[8pt]
&=
\displaystyle \sum_{i=1}^r \sum_{j=1}^r a_{ij}{\rm Tr}
(P_{\alpha_i}^THP_{\alpha_i}){\rm Tr}
(P_{\alpha_j}^THP_{\alpha_j})\\[8pt]
&\quad \, +\displaystyle  \sum_{i=1}^r \sum_{j=1}^r \delta_{ij}\displaystyle \sum_{i'\in \alpha_i,j'\in \alpha_j}(p_{i'}^THp_{j'})^2.
\end{array}
$$
The proof is completed. \hfill $\Box$
\begin{corollary}\label{coro:DFz}
Let $w_1(z):=u_{k_1(z)}(z),\ldots, w_{r(z)}(z):=u_{k_r(z)}(z)$ be $r(z)$ distinct values of $m$ values  of $u(z)$, and $X(z)$ has the spectral decomposition $X(z)=P{\rm Diag}u(z) P^T$ with $u(z)=\lambda(X(z))$,
namely
$$
u_1(z)=\cdots=u_{k_1(z)}(z) > u_{k_1(z)+1}(z)=\cdots=u_{k_2(z)}(z) > u_{k_2(z)+1}(z)\cdots u_{k_r(z)}(z),
$$
where $k_0=0, k_r(z)=m$. Denote
$$
\alpha_1(z)=\{1,\ldots, k_1(z)\}, \alpha_2(z)=\{k_1(z)+1,\ldots, k_2(z)\}, \ldots, \alpha_{r(z)}(z)=\{k_{r-1}(z)+1,\ldots, k_{r(z)}\}.
$$
Let $\psi: \mathbb R^m \rightarrow \mathbb R$ be a symmetric function, twice differentiable at the point $w \in \mathbb R^m$,
and $P$ be a permutation matrix such that $Pw=w$. Then
\begin{itemize}
\item[(1)]there are real numbers $b_1(z),\ldots, b_{r(z)}(z)$ and a symmetric matrix $(a_{ij}(z))_{i,j=1}^{r(z)}$ such that
\begin{equation}\label{2ndpsimuz}
\begin{array}{l}
\nabla^2 \psi (u(z))=\left[
\begin{array}{ccc}
b_1(z)I_{|\alpha_1(z)|} & & \\[4pt]

 & \ddots  &\\[4pt]
& & b_{r(z)}(z)I_{|\alpha_{r(z)}(z)|}
\end{array}
\right]\\[18pt]
+\left [
\begin{array}{ccc}
a_{11}(z)\textbf{1}_{|\alpha_1(z)|}\textbf{1}_{|\alpha_1(z)|}^T &\cdots & a_{1r(z)}(z)\textbf{1}_{|\alpha_1(z)|}\textbf{1}_{|\alpha_{r(z)}(z)|}^T\\[10pt]
 \vdots &\vdots& \vdots\\[10pt]
a_{{r(z)}1}(z)\textbf{1}_{|\alpha_{r(z)}(z)|}\textbf{1}_{|\alpha_1(z)|}^T & \cdots & a_{r(z)r(z)}(z)\textbf{1}_{|\alpha_{r(z)}(z)|}\textbf{1}_{|\alpha_{r(z)}(z)|}^T
\end{array}
\right].
\end{array}
\end{equation}
\item[(2)]there are real numbers $c_1(z),\ldots, c_{r(z)}(z)$ such that
$$
\nabla \psi (u(z))=\left[
\begin{array}{c}
c_1(z)\textbf{1}_{|\alpha_1(z)|}\\[2pt]
\vdots\\[2pt]
c_r(z)\textbf{1}_{|\alpha_{r(z)}(z)|}
\end{array}
\right],
$$
then for any $H \in \mathbb S^m$,
\begin{equation}\label{eq:2ndH}
\begin{array}{l}
\nabla^2 (\psi\circ\lambda)(X(z))[H]=\\[10pt]
P\left[
\begin{array}{ccc}
\left(\displaystyle \sum_{j=1}^{r(z)}a_{1j}(z) {\rm Tr}\left(P_{\alpha_j(z)}^THP_{\alpha_j(z)}\right)\right)I_{|\alpha_1(z)|}  & &\\[4pt]
&  \ddots  &\\[4pt]
 & & \left(\displaystyle \sum_{j=1}^{r(z)}a_{r(z)j}(z) {\rm Tr}\left(P_{\alpha_j(z)}^THP_{\alpha_j(z)}\right)\right) I_{|\alpha_{r(z)}(z)|}
\end{array}
\right]P^T\\[12pt]
+P\left[
\begin{array}{ccc}
\delta_{11}(z)P_{\alpha_1(z)}^THP_{\alpha_1(z)}  & \cdots & \delta_{1r}(z)P_{\alpha_1(z)}^THP_{\alpha_{r(z)}(z)}\\[4pt]
\vdots &  \ddots  & \vdots\\[4pt]
\delta_{r(z)1}P_{\alpha_{r(z)}(z)}^THP_{\alpha_1(z)}  & \cdots & \delta_{r(z)r(z)}P_{\alpha_{r(z)}(z)}^THP_{\alpha_{r(z)}(z)}
\end{array}
\right]P^T
\end{array}
\end{equation}
and
\begin{equation}\label{eq:2ndDSF}
\begin{array}{l}
\nabla^2 (\psi\circ\lambda)(A)[H,H]\\[6pt]
 =\displaystyle \sum_{i=1}^{r(z)} \sum_{j=1}^{r(z)} a_{ij}(z){\rm Tr}
(P_{\alpha_i(z)}^THP_{\alpha_i(z)}){\rm Tr}
(P_{\alpha_j(z)}^THP_{\alpha_j(z)})\\[12pt]
\quad \quad +\displaystyle \sum_{i=1}^{r(z)} \sum_{j=1}^{r(z)} \delta_{ij}(z)\displaystyle \sum_{i'\in \alpha_i(z),j'\in \alpha_j(z)}(p_{i'}^THp_{j'})^2,
\end{array}
\end{equation}
\end{itemize}
where
$$
\left\{
\begin{array}{l}
\delta_{ii}(z)=b_i(z), i=1,\ldots, r(z),\\[8pt]
\delta_{ij}(z)=\displaystyle \frac{c_i(z)-c_j(z)}{w_i(z)-w_j(z)},i,j=1,\ldots r(z), i \ne j.
\end{array}
\right.
$$
\end{corollary}

\begin{proposition}\label{prop2nd-der}
Under the setting of Proposition \ref{prop:DF}, for any $H\in \mathbb S^m$, one has
$Y=G(z)^{-1}H$ can be characterized by
\begin{equation}\label{eq:Yexp}
\begin{array}{l}
P_{\alpha_i(z)}^TYP_{\alpha_j(z)}=\displaystyle \frac{1}{1+\mu(z)\delta_{ij}(z)}P_{\alpha_i(z)}^THP_{\alpha_j(z)},i,j=1,\ldots, r(z), i\ne j,
\end{array}
\end{equation}
for $k=1,\ldots, r(z)$, $i,j \in \alpha_k$, $i\ne j$,
\begin{equation}\label{eq:Yexpij}
\begin{array}{l}
P_i^TYP_j=\displaystyle \frac{1}{1+\mu(z)\delta_{kk}(z)}P_i^THP_j,
\end{array}
\end{equation}
and
\begin{equation}\label{eq:Yexpa}
{\rm diag}(P^TYP)=[I_m+\mu(z) \nabla^2 \psi (u(z))]^{-1}
{\rm diag}(P^THP).
\end{equation}
\end{proposition}
{\bf Proof}. For $Y=G(z)^{-1}H$, we have
$$
\begin{array}{l}
H=G(z)Y =Y+\mu(z)\nabla^2(\varphi_p \circ \lambda)(X(z))[Y]\\[12pt]
=Y+\mu(z)
P\left[
\begin{array}{ccc}
\left(\displaystyle \sum_{j=1}^{r(z)}a_{1j}(z) {\rm Tr}\left(P_{\alpha_j(z)}^TYP_{\alpha_j(z)}\right)\right)I_{|\alpha_1(z)|}  & &\\[4pt]
&  \ddots  &\\[4pt]
 & & \left(\displaystyle \sum_{j=1}^{r(z)}a_{rj}(z) {\rm Tr}\left(P_{\alpha_j(z)}^TYP_{\alpha_j(z)}\right)\right) I_{|\alpha_{r(z)}(z)|}
\end{array}
\right]P^T\\[12pt]
+\mu(z)P\left[
\begin{array}{ccc}
\delta_{11}(z)P_{\alpha_1(z)}^TYP_{\alpha_1(z)}  & \cdots & \delta_{1r}(z)P_{\alpha_1(z)}^TYP_{\alpha_{r(z)}(z)}\\[4pt]
\vdots &  \ddots  & \vdots\\[4pt]
\delta_{r(z)1}P_{\alpha_{r(z)}(z)}^TYP_{\alpha_1(z)}  & \cdots & \delta_{r(z)r(z)}P_{\alpha_{r(z)}(z)}^TYP_{\alpha_{r(z)}(z)}
\end{array}
\right]P^T.
\end{array}
$$
This implies, for $i,j=1,\ldots, r(z), i\ne j$,
$$
P_{\alpha_i(z)}^TYP_{\alpha_j(z)}=\displaystyle \frac{1}{1+\mu(z)\delta_{ij}(z)}P_{\alpha_i(z)}^THP_{\alpha_j(z)}.
$$
And for $k=1,\ldots, r(z)$, we have
$$
\begin{array}{ll}
P_{\alpha_k(z)}^THP_{\alpha_k(z)} &=(1+\mu(z)\delta_{kk}(z))P_{\alpha_k(z)}^TYP_{\alpha_k(z)}+I_{\alpha_k(z)}\mu(z)
\left(\displaystyle \sum_{j=1}^{r(z)}a_{kj}(z) {\rm Tr}\left(P_{\alpha_j(z)}^TYP_{\alpha_j(z)}\right)\right)\\[12pt]
&=(1+\mu(z)\delta_{kk}(z))P_{\alpha_k(z)}^TYP_{\alpha_k(z)}+\mu(z)
\left(\displaystyle \sum_{j=1}^{r(z)}a_{kj}(z) I_{\alpha_k(z)}I_{\alpha_j(z)}^*\left(P_{\alpha_j(z)}^TYP_{\alpha_j(z)}\right)\right)
\end{array}
$$
which implies for $k=1,\ldots, r(z)$, $i,j \in \alpha_k(z)$, $i\ne j$,
$$
P_i^TYP_j=\displaystyle \frac{1}{1+\mu(z)\delta_{kk}(z)}P_i^THP_j.
$$
For $k=1,\ldots, r(z)$,
\begin{equation}\label{eq:reqs}
\begin{array}{ll}
{\rm diag}[P_{\alpha_k(z)}^THP_{\alpha_k(z)}] &=(1+\mu(z)\delta_{kk}(z)){\rm diag}[P_{\alpha_k(z)}^TYP_{\alpha_k(z)}]\\[12pt]
& \quad +\mu(z)\textbf{1}_{|\alpha_k(z)|}
\left(\displaystyle \sum_{j=1}^{r(z)} a_{ij}(z) \textbf{1}_{|\alpha_j(z)|}^T {\rm diag}\left(P_{\alpha_j(z)}^TYP_{\alpha_j(z)}\right)\right).
\end{array}
\end{equation}
Define
$$
\left(\begin{array}{c}
\widehat y_1\\[4pt]
\widehat y_2\\[4pt]
\vdots\\[4pt]
\widehat y_{r(z)}
\end{array}
\right)=
\left(\begin{array}{c}
{\rm diag}\left(P_{\alpha_1(z)}^TYP_{\alpha_1(z)}\right)\\[4pt]
{\rm diag}\left(P_{\alpha_2(z)}^TYP_{\alpha_2(z)}\right)\\[4pt]
\vdots\\[4pt]
{\rm diag}\left(P_{\alpha_{r(z)}(z)}^TYP_{\alpha_{r(z)}(z)}\right)
\end{array}
\right)\mbox{ and }\left(\begin{array}{c}
\widehat h_1\\[4pt]
\widehat h_2\\[4pt]
\vdots\\[4pt]
\widehat h_{r(z)}
\end{array}
\right)=
\left(\begin{array}{c}
{\rm diag}\left(P_{\alpha_1(z)}^THP_{\alpha_1(z)}\right)\\[4pt]
{\rm diag}\left(P_{\alpha_2(z)}^THP_{\alpha_2(z)}\right)\\[4pt]
\vdots\\[4pt]
{\rm diag}\left(P_{\alpha_{r(z)}(z)}^THP_{\alpha_{r(z)}(z)}\right)
\end{array}
\right).
$$
Then the solution $(\widehat y_1,\ldots, \widehat y_r)$ of (\ref{eq:reqs}) is the solution to the following system of linear equations:
$$
\begin{array}{l}
\left(\begin{array}{c}
\widehat h_1\\[4pt]
\vdots\\[4pt]
\widehat h_r
\end{array}
\right)=M(z)
 \left(
 \begin{array}{c}
\widehat y_1\\[4pt]
\vdots\\[4pt]
\widehat y_{r(z)}
\end{array}
\right)
\end{array}
$$
where
$$
\begin{array}{ll}
M(z)=&
\left\{
\left( \begin{array}{ccc}
\widehat \delta_{11}(z)I_{|\alpha_1(z)|}
 & \cdots & 0\\[12pt]
\vdots & \vdots & \vdots\\[12pt]
0 & \cdots &\widehat \delta_{rr}(z))I_{|\alpha_{r(z)}(z)|}
 \end{array}
 \right)\right.\\[12pt]
&+\left.\left(\begin{array}{ccc}
\widehat a_{11}(z)\textbf{1}_{|\alpha_1(z)|}\textbf{1}_{|\alpha_1 (z)|}^T
 & \cdots & \widehat a_{1r(z)}(z)\textbf{1}_{|\alpha_1(z)|}\textbf{1}_{|\alpha_{r(z)}(z)|}^T\\[12pt]
\vdots & \vdots & \vdots \\[12pt]
\widehat a_{r1}(z)\textbf{1}_{|\alpha_{r(z)}(z)|}\textbf{1}_{|\alpha_1(z)|}^T
 & \cdots &\widehat a_{rr}(z))\textbf{1}_{|\alpha_{r(z)}(z)|}\textbf{1}_{|\alpha_{r(z)}(z)|}^T
 \end{array}
 \right)\right\}
 \end{array}
$$
with
$$
\widehat \delta_{ii}(z)=1+\mu(z)\delta_{ii}(z), i=1,\ldots, r(z), \, \widehat a_{ij}(z)=\mu(z)a_{ij}(z), i,j=1,\ldots, r(z).
$$
It is not difficult to verify that
\begin{equation}\label{eq:Mzdef}
M(z)=I_m+\mu(z) \nabla^2 \psi (u(z)).
\end{equation}
Therefore, the formulas (\ref{eq:Yexp}), (\ref{eq:Yexpij}) and (\ref{eq:Yexpa}) hold. \hfill $\Box$\\
Define
$$
\pi_{ij}(z)=(1+\mu(z)\delta_{ij}(z))^{-1}, i,j=1,\ldots, r(z)
$$
and
\begin{equation}\label{2ndpsimuB}
\begin{array}{ll}
{\cal B}(z)=&-\left[
\begin{array}{ccc}
\pi_{11}(z)I_{|\alpha_1(z)|}  & &\\[4pt]
 & \ddots  &\\[4pt]
 & & \pi_{r(z)r(z)}(z)I_{|\alpha_r(z)|}
\end{array}
\right]\\[18pt]
& +\left [
\begin{array}{cccc}
\pi_{11}(z)\textbf{1}_{|\alpha_1(z)|}\textbf{1}_{|\alpha_1(z)|}^T &\cdots & \pi_{1r(z)}(z)\textbf{1}_{|\alpha_1(z)|}\textbf{1}_{|\alpha_{r(z)}(z)|}^T\\[10pt]
 \vdots &\vdots& \vdots\\[10pt]
\pi_{r(z)1}(z)\textbf{1}_{|\alpha_{r(z)}(z)|}\textbf{1}_{|\alpha_1(z)|}^T & \cdots & \pi_{r(z)r(z)}(z)\textbf{1}_{|\alpha_{r(z)}(z)|}\textbf{1}_{|\alpha_{r(z)}(z)|}^T
\end{array}
\right].
\end{array}
\end{equation}
Then, for any $H \in \mathbb S^m$,
\begin{equation}\label{eq:gnyH}
G^{-1}(z)H=P\left[ {\rm Diag}\left(M(z)^{-1}{\rm diag}(P^THP)\right)+{\cal B}(z)\circ  P^THP\right]P^T.
\end{equation}

\begin{proposition}\label{prop:important}
Let $\psi: \mathbb R^m \rightarrow \mathbb R_+$ be a twice differentiable symmetric gauge function. Let $z=(A,s)\in {\cal H}$ be a given point with $A=P{\rm Diag}(\lambda (A))P^T$. Let $X(z)=P {\rm Diag}(u(z))P^T$ and $(u(z),\mu(z))\in \mathbb R^m \times \mathbb R_{++}$  is a solution of the system of equations:
$$
F(u,\mu;z)=0.
$$
Then for
  $$\beta(z)=(1+\nabla (\psi\circ \lambda)(X(z))^*G(z)^{-1}\nabla (\psi\circ \lambda)(X(z)))^{-1},\, G(z)={\cal I}+\mu(z) \nabla^2 (\psi\circ \lambda)(X(z)),$$
we have
\begin{equation}\label{eq:Ginvnab}
G(z)^{-1}\nabla (\psi \circ \lambda) (X(z))=P {\rm Diag}\left[(I_m+\mu(z)\nabla^2\psi(u(z)))^{-1}\nabla \psi(u(z))\right]P^T
\end{equation}
and
\begin{equation}\label{eq:betaz}
\beta(z)=1+\nabla \psi(u(z))^T(I_m+\mu(z)\nabla^2\psi(u(z)))^{-1}\nabla \psi(u(z)).
\end{equation}
\end{proposition}
\begin{remark}\label{remark}
Proposition \ref{prop:important} gives elegant expressions for $G(z)^{-1}\nabla (\psi \circ \lambda) (X(z))$ and $\beta(z)$, which appear in Proposition \ref{mpropPiD}, Proposition \ref{mpropqPiDq}, Theorem \ref{mthPiDirction} and Theorem \ref{mthBsubdiff}.
\end{remark}
\section{Variational Analysis on Schatten $p$-norm Cone}\label{Section-Kn}
\setcounter{equation}{0}
In this section, we consider a special norm of $\mathbb S^m$,  the Schatten $p$-norm, induced by
$p$-norm of $\mathbb R^m$ defined by
$$
\varphi_p(u)=\left(\displaystyle \sum_{j=1}^m |u_j|^p\right)^{1/p},
$$
for $p\in (1,+\infty)$ is a real number. When $p=1$,
$$
\varphi_1(u)=\displaystyle \sum_{j=1}^m |u_j|
$$
and when $p=+\infty$,
$$
\varphi_{\infty}(u)=\displaystyle \max_{1\leq j \leq m} |u_j|.
$$
We use  $\K_p$ to denote the Schatten $p$-norm cone for $p \in [0,\infty]$, namely
\begin{equation}\label{eq:spn}
\K_p=\Big \{(A,s)\in {\cal H}: (\varphi_p\circ \lambda)(A)\leq s\Big\}.
\end{equation}
For $r \in (1,\infty)$, the conjugate number of $r$, denoted by $r_*$, is defined through
$$
\displaystyle \frac{1}{r}+\displaystyle \frac{1}{r_*}=1.
$$
Obviously we have that, if $N(\cdot)=(\varphi_r\circ \lambda)(\cdot)$, the dual norm is $N_*(\cdot)=(\varphi_{r_*}\circ \lambda)(\cdot)$.
Now we introduce some notations which will be used in the following analysis. Denote for $w \in \Re^m$
$$
\begin{array}{l}
|w|=(|w_1|,\ldots, |w_m|)^T\\[4pt]
|w|^r=(|w_1|^r,\ldots, |w_m|^r)^T\mbox{ for } r \in \mathbb R\\[4pt]
{\rm sgn} (z)=({\rm sgn}(z_1),\ldots, {\rm sgn}(z_m))^T\\[4pt]
W={\rm diag}(w_1, \ldots, w_m)\\[4pt]
|W|={\rm diag}(|w_1|, \ldots, |w_m|)\\[4pt]
|W|^r={\rm diag}(|w_1|^r, \ldots, |w_m|^r)\mbox{ for } r \in \mathbb R\\[4pt]
a\circ b=(a_1b_1,\ldots a_mb_m)^T \mbox{ for } a,b \in \mathbb R^m\\[4pt]
\textbf{1}_m=(1,\ldots, 1)^T \in \mathbb R^m
\end{array}
$$
It is easy to see that, for $p=1$ or $p=\infty$, the corresponding cones, ${\cal K}_1$ and ${\cal K}_{\infty}$, are nuclear norm cone and spectral norm cone, respectively.  In Subsection \ref{Sec3.1}, we will discuss the variational geometry of  ${\cal K}_1$ and ${\cal K}_{\infty}$. In Subsections \ref{Sec3.2}-\ref{Sec3.4}, we will study variational properties of matrix norm cone ${\cal K}_p$ for $p \in (1,\infty)$.
\subsection{Variational Analysis of ${\cal K}_1$ and ${\cal K}_{\infty}$ }\label{Sec3.1}
As we mentioned,Ding (2017) \cite{Ding2017} studied  some variational properties of the spectral,
and nuclear matrix norm.  However, for completeness of variational analysis of Schatten $p$-norm cones,
we present tangent cones, normal cones, outer second-order tangent sets of ${\cal K}_1$ and ${\cal K}_{\infty}$,
as well as expressions of projection operators on ${\cal K}_1$ and ${\cal K}_{\infty}$. For directional derivatives of projection operators on ${\cal K}_1$ and ${\cal K}_{\infty}$, see the related results in  \cite{Ding2017}.
\subsubsection{Variational geometry of ${\cal K}_1$}
Let $\phi_1: \mathbb S^{m}\times \mathbb R$ be defined by
$$
\phi_1 (A,s)=(\varphi_1\circ \lambda)(A)-s.
$$
 Then $\K_1$  can be expressed as the following level set
 $$
 \K_1=\{(A,s)\in \mathbb S^{m}\times \mathbb R: \phi_1(A,s)\leq 0\}.
 $$
  For $z=(A,s)\ne (0,0)$, $\phi_1(A,s)=0$, define
 \begin{equation}\label{eq:index1}
 \alpha=\{i: \lambda_i(A)>0\},\, \beta=\{i: \lambda_i(A)=0\},\, \gamma=\{i:\lambda_i(A) <0\}.
 \end{equation}
 Then the directional derivative of $\phi_1$ at $z$ is
 $$
 \begin{array}{ll}
 \phi_1'(z;d_z)&=\varphi'(\lambda(A);\lambda'(A;d_A))-d_s\\[4pt]
 &=\displaystyle \sum_{i\in \alpha} \lambda'_i(A;d_A)-\displaystyle \sum_{i\in \gamma} \lambda'_i(A;d_A)+
 \displaystyle \sum_{i\in \beta} |\lambda'_i(A;d_A)|-d_s\\[4pt]
 &= {\rm Tr}\,(P_{\alpha}^Td_AP_{\alpha})-{\rm Tr}\,(P_{\gamma}^Td_AP_{\gamma})+
 (\varphi^{|\beta|}_1\circ \lambda)(P_{\beta}^Td_AP_{\beta})-d_s,
 \end{array}
 $$
 where
 $$\varphi^{|\beta|}_1(v)=\displaystyle\sum_{i=1}^{|\beta|}|v_i| \mbox{ for } v \in \mathbb R^{|\beta|}.
 $$
 For convenience of discussions, without loss of generality, we assume
 $$
 \alpha=\{1,\ldots, |\alpha|\}, \beta=\{|\alpha|+1,\ldots, |\alpha|+|\beta|\}, \, \gamma=\{ |\alpha|+|\beta|,\ldots, m\}
 $$
 and
 $$
 \alpha=\alpha_1\cup \cdots \cup \alpha_{r_1-1}, \beta=\alpha_{r_1}, \gamma=\alpha_{r_1+1}\cup \cdots \cup \alpha_{r}.
 $$
 When $\phi_1(z)=0$ and $\phi_1'(z,d_z)=0$, denote
 \begin{equation}\label{eq:index2}
 \begin{array}{l}
  \beta_+=\{|\alpha|+i \in \beta: \lambda_i (P_{\beta}^Td_AP_{\beta})>0\},\\[3pt]
  \beta_0=\{|\alpha|+i \in \beta: \lambda_i (P_{\beta}^Td_AP_{\beta})=0\},\\[3pt]
  \beta_-=\{|\alpha|+i \in \beta: \lambda_i (P_{\beta}^Td_AP_{\beta})<0\}.
    \end{array}
 \end{equation}
 Let ${\cal Q}\in {\cal O}^{|\beta|}$ with $Q=(Q_{ \beta_+} \, Q_{ \beta_0}\,Q_{ \beta_-})$ such that
 $$
 P_{\beta}^Td_AP_{\beta}=Q{\rm Diag}\,(\lambda(P_{\beta}^Td_AP_{\beta}))Q^T.
 $$
  In this case, the second-order directional derivative of $\phi_1$ at $z$ along $(d_z,\xi_z)$ is
 $$
 \begin{array}{l}
  \phi_1''(z;d_z,\xi_z)=\varphi''(\lambda(A);\lambda'(A;d_A),\lambda''(A;d_A,\xi_A))-\xi_s\\[8pt]
  =\displaystyle \sum_{k=1}^{r_1-1}{\rm Tr}\,\left(P_{\alpha_k}^T[\xi_A-2d_A(A-\mu_kI_m)^{\dagger}d_A]P_{\alpha_k}\right)
  +{\rm Tr}\,\left(Q_{ \beta_+}^TP_{\beta}^T[\xi_A-2d_AA^{\dagger}d_A]P_{\beta}Q_{ \beta_+}\right)  \\[16pt]
  -\displaystyle \sum_{k=r_1+1}^{r}{\rm Tr}\,\left(P_{\alpha_k}^T[\xi_A-2d_A(A-\mu_kI_m)^{\dagger}d_A]P_{\alpha_k}\right)
  -{\rm Tr}\,\left(Q_{ \beta_-}^TP_{\beta}^T[\xi_A-2d_AA^{\dagger}d_A]P_{\beta}Q_{ \beta_-}\right)  \\[16pt]
   + (\varphi^{|\beta_0|}_1 \circ \lambda)
 \left(Q_{ \beta_0}^TP_{\beta}^T[\xi_A-2d_AA^{\dagger}d_A]P_{\beta}Q_{ \beta_0}\right) -\xi_s.
 \end{array}
 $$
 We may use Proposition 2.61 and Proposition 3.30 of \cite{BS00} to derive the tangent cone and the second-order tangent set of $\K_1$, respectively.
 \begin{proposition}\label{VAl1norm}
 The tangent cone and the second-order tangent set can be characterized by the following assertions:
  \begin{itemize}
 \item[{\rm (i)}] For any $(A,s)\in \K^{m+1}_1$,
$$
{\cal T}_{\K_1}(A,s)=\left
\{
\begin{array}{ll}
\mathbb S^{m}\times \mathbb R &  (\varphi_1\circ\lambda)(A) < s,\\[2pt]
\K_1 & (A,s)=(0,0),\\[2pt]
\left\{(d_A,d_s): \begin{array}{l}
{\rm Tr}\,(P_{\alpha}^Td_AP_{\alpha})-{\rm Tr}\,(P_{\gamma}^Td_AP_{\gamma})\\[6pt]
+ (\varphi^{|\beta|}_1\circ \lambda)(P_{\beta}^Td_AP_{\beta})\leq d_s \end{array}\right\}
 & (\varphi_1\circ\lambda)(A)=s>0.
\end{array}
\right.
$$
\item[{\rm (ii)}]
Let $z \in \K_1$ and $d \in {\cal T}_{ \K_1}(z)$ where $z=(A,s)$ and $d=(d_w,d_s)$. Then
\begin{equation}\label{eq:l1normCT}
{\cal T}^2_{ \K_1}(z,d)=\left\{
\begin{array}{ll}
\mathbb S^{m}\times \mathbb R & d \in \mbox{int }{\cal T}_{ \K_1}(z),\\[6pt]
{\cal T}_{ \K_1}(z) & z=0,\\[6pt]
\left\{(\xi_A,\xi_s): \begin{array}{l}\displaystyle \sum_{k=1}^{r_1-1}{\rm Tr}\,\left(P_{\alpha_k}^T[\xi_A-2d_A(A-\mu_kI_m)^{\dagger}d_A]P_{\alpha_k}\right)\\[6pt]
  +{\rm Tr}\,\left(Q_{ \beta_+}^TP_{\beta}^T[\xi_A-2d_AA^{\dagger}d_A]P_{\beta}Q_{ \beta_+}\right)  \\[6pt]
  -\displaystyle \sum_{k=r_1+1}^{r}{\rm Tr}\,\left(P_{\alpha_k}^T[\xi_A-2d_A(A-\mu_kI_m)^{\dagger}d_A]P_{\alpha_k}\right)\\[16pt]
  -{\rm Tr}\,\left(Q_{ \beta_-}^TP_{\beta}^T[\xi_A-2d_AA^{\dagger}d_A]P_{\beta}Q_{ \beta_-}\right)  \\[16pt]
   + (\varphi^{|\beta_0|}_1 \circ \lambda)
 \left(Q_{ \beta_0}^TP_{\beta}^T[\xi_A-2d_AA^{\dagger}d_A]P_{\beta}Q_{ \beta_0}\right) -\xi_s\leq0
  \end{array}
  \right\} & \mbox{otherwise}.
\end{array}
\right.
\end{equation}
\end{itemize}
\end{proposition}
 \subsubsection{Variational geometry of ${\cal K}_{\infty}$}
 Let $\phi_{\infty}: \mathbb S^{m}\times \mathbb R\rightarrow \mathbb R$ be defined by
$$
\phi_{\infty} (A,s)=(\varphi_{\infty}\circ \lambda)(A)-s.
$$
 Then $\K_{\infty}$ can be expressed as the following level set
 $$
  \K_{\infty}=\{(A,s)\in \mathbb S^{m}\times \mathbb R: \phi_{\infty}(A,s)\leq 0\}.
 $$
  For $z=(A,s)\ne (0,0)$, $\phi_{\infty}(A,s)=0$, define
 \begin{equation}\label{eq:index1a}
 \alpha_+=\{i: \lambda_i(A)=s\},\, \alpha_-=\{i:\lambda_i(A)=-s\},\, \gamma=\{i: |\lambda_i(A)| <s\}.
 \end{equation}
 Then the directional derivative of $\phi_{\infty}$ at $z$ is
 $$
 \begin{array}{ll}
 \phi_{\infty}'(z;d_z)&=\displaystyle \max \left[\displaystyle\max_{i\in \alpha_+} \lambda_i'(A;d_{A}),\,\,\displaystyle\max_{i\in \alpha_-}[- \lambda_i'(A;d_{A})]\right]-d_s\\[6pt]
 &=\displaystyle \max \left[\lambda_{\max}(P_{\alpha_+}^Td_AP_{\alpha_+}),\,\,-\lambda_{\min}(P_{\alpha_-}^Td_AP_{\alpha_-})\right]-d_s\\[6pt]
 \end{array}
 $$
 When $\phi_{\infty}(z)=0$ and $\phi_{\infty}'(z,d_z)=0$, denote
 \begin{equation}\label{eq:index2a}
 \alpha^*_+=\{i \in \alpha_+: \lambda_i'(A;d_{A})=d_s\},\,\alpha^*_-=\{i \in \alpha_-: \lambda_i'(A;d_{A})=-d_s\}.
 \end{equation}
 Let $Q^+\in {\cal O}^{|\alpha_+|}$ and $Q^-\in {\cal O}^{|\alpha_-|}$ be such that
 $$
 P_{\alpha_+}^Td_AP_{\alpha_+}=Q^+{\rm Diag}\,\lambda(P_{\alpha_+}^Td_AP_{\alpha_+})Q^{+T},\,
  P_{\alpha_-}^Td_AP_{\alpha_-}=Q^-{\rm Diag}\,\lambda(P_{\alpha_-}^Td_AP_{\alpha_-})Q^{-T}.
 $$

 In this case, the second-order directional derivative of $\phi_{\infty}$ at $z$ along $(d_z,\xi_z)$ is
 $$
 \begin{array}{ll}
  \phi_{\infty}''(z;d_z,\xi_z)&=\displaystyle \max \left[\displaystyle\max_{i\in \alpha^*_+} \lambda''_i(A;d_A,\xi_A),\,\,\displaystyle\max_{i\in \alpha^*_-}[-  \lambda''_i(A;d_A,\xi_A)]\right]-\xi_s\\[8pt]
  &=\displaystyle \max \left[\lambda_{\max}\left((Q^+_{\alpha^*_+})^T P_{\alpha_+}^T[\xi_A-2d_A(A-sI_m)^{\dagger}d_A] P_{\alpha_+}Q^+_{\alpha^*_+})\right),\right.\\[8pt]
  & \left.\quad \quad \quad \, -\lambda_{\min}\left((Q^-_{\alpha^*_-})^T P_{\alpha_-}^T[\xi_A-2d_A(A+sI_m)^{\dagger}d_A] P_{\alpha_-}Q^-_{\alpha^*_-})\right)\right]-\xi_s.
  \end{array}
 $$

  Again, we may use Proposition 2.61 and Proposition 3.30 of \cite{BS00} to derive the tangent cone and the second-order tangent set of $\K_{\infty}$, respectively.
 \begin{proposition}\label{VAlinfnorm}
 The tangent cone and the second-order tangent set can be characterized by the following assertions:
  \begin{itemize}
 \item[{\rm (i)}] For any $(A,s)\in \K_{\infty}$,
$$
{\cal T}_{\K_{\infty}}(A,s)=\left
\{
\begin{array}{ll}
\mathbb S^{m}\times \mathbb R &  (\varphi_{\infty}\circ\lambda)(A) < s,\\[2pt]
\K_{\infty} & (A,s)=(0,0),\\[2pt]
\left\{(d_w,d_s)\in \mathbb S^{m}\times \mathbb R:
\begin{array}{l}
P_{\alpha_+}^Td_AP_{\alpha_+} \preceq d_s I_{|\alpha_+|}\\[4pt]
P_{\alpha_-}^Td_AP_{\alpha_-} \succeq -d_s I_{|\alpha_-|}
\end{array}
\right\} & (\varphi_{\infty}\circ\lambda)(A)=s>0.
\end{array}
\right.
$$
\item[{\rm (ii)}]
Let $z \in \K_{\infty}$ and $d \in {\cal T}_{ \K_{\infty}}(z)$ where $z=(A,s)$ and $d=(d_w,d_s)$. Then
\begin{equation}\label{eq:linfnormCT}
{\cal T}^2_{ \K_{\infty}}(z,d)=\left\{
\begin{array}{ll}
\mathbb S^{m}\times \mathbb R & d \in \mbox{int }{\cal T}_{ \K_{\infty}}(z),\\[6pt]
{\cal T}_{ \K_{\infty}}(z) & z=0,\\[6pt]
\left\{(\xi_w,\xi_s): \begin{array}{l}
(Q^+_{\alpha^*_+})^T P_{\alpha_+}^T[\xi_A-2d_A(A-sI_m)^{\dagger}d_A] P_{\alpha_+}Q^+_{\alpha^*_+})\preceq \xi_s I_{|\alpha^*_+|}\\[6pt]
(Q^-_{\alpha^*_-})^T P_{\alpha_-}^T[\xi_A-2d_A(A+sI_m)^{\dagger}d_A] P_{\alpha_-}Q^-_{\alpha^*_-}) \succeq
-\xi_s I_{|\alpha^*_-|}
\end{array}
\right\} & \mbox{otherwise}.
\end{array}
\right.
\end{equation}
\end{itemize}
\end{proposition}
 \subsubsection{Projections over $\K_1$ and  $\K_{\infty}$}
Define
$$
{\cal P}_{\mu}(y)=[|y|-\mu \textbf{1}_m]_+ \circ {\rm sgn}\,(y)
$$
for $\mu \in \mathbb R$ and $y \in \mathbb R^m$.

It follows from Example 6.38 of \cite{Beck2017} that
\begin{lemma}\label{lem:projl1}
Let $u=(a,s)\in \mathbb{R}^{m}\times \mathbb{R}$ be a given point, then
\begin{equation}\label{eq:kfap}
\Pi_{{\rm epi}\,\varphi_1}(a,s)=\left
\{
\begin{array}{ll}
(a,s) & \mbox{if } (a,s) \in {\rm epi}\,\varphi_1,\\[2pt]
(0,0) & \mbox{if } (a,s) \in [{\rm epi}\,\varphi_1]^{\circ},\\[2pt]
({\cal P}_{\mu (u)}(a), \mu(u)+s) &  \mbox{otherwise},
\end{array}
\right.
\end{equation}
where $\mu(u)\in \mathbb R_{++}$ is any positive root of the nonincreasing function:
$$
\theta (\mu)=\|{\cal P}_{\mu}(a)\|_1 -\mu -s.
$$
\end{lemma}
Noting that
$$
\theta (\mu)=\displaystyle \sum_{i=1}^m {[|a_i|-\mu]}_+ - \mu -s,
$$
we have, for  $|s| < \|a\|_1$, that
$$
\left\{
\begin{array}{ll}
\theta (0)=\|a\|_1-s>0, \, \theta (\|a\|_1)=-\|a\|_1-s<0 & \mbox{if } s \geq 0,\\[6pt]
\theta (0)=\|a\|_1-s>0, \, \theta (\mu_0)=-\mu_0-s\leq -1 \mbox{ for }\mu_0=
\max\{\|a\|_1,-s\}+1  & \mbox{if } s < 0.
\end{array}
\right.
$$
From the monotonicity of $\theta$, we have that
$$
\mu(u) \in \left\{
\begin{array}{ll}
(0, \|a\|_1) & \mbox{ if } s\geq0,\\[5pt]
(0, \mu_0) & \mbox{ if } s<0.
\end{array}
\right.
$$
From the identity
$$
u=\Pi_{{\rm epi}\, \varphi_{\infty}}(u)+\Pi_{[{\rm epi}\, \varphi_{\infty}]^{\circ}}(u)=\Pi_{{\rm epi}\, \varphi_{\infty}}(u)+\Pi_{-{\rm epi}\, \varphi_1}(u),
$$
and $\Pi_{-{\rm epi}\, \varphi_1}(u)=\Pi_{{\rm epi}\, \varphi_1}(-u)$, we have
\begin{equation}\label{eq:linfp}
\Pi_{{\rm epi}\, \varphi_{\infty}}(u)=u-\Pi_{{\rm epi}\, \varphi_1}(-u).
\end{equation}
From (\ref{eq:linfp}), we obtain the formula for the projection onto ${\rm epi}\, \varphi_{\infty}$.
\begin{lemma}\label{lem:projlinf}
Let $u=(a,s)\in \mathbb{R}^{m}\times \mathbb{R}$ be a given point, then
\begin{equation}\label{eq:kfap}
\Pi_{{\rm epi}\, \varphi_{\infty}}(a,s)=\left
\{
\begin{array}{ll}
(a,s) & \mbox{if } (a,s) \in {\rm epi}\, \varphi_{\infty},\\[2pt]
(0,0) & \mbox{if } (a,s) \in -{\rm epi}\, \varphi_1,\\[2pt]
(a,s)-({\cal P}_{\mu (-u)}(-a), \lambda(-u)-s) &  \mbox{otherwise},
\end{array}
\right.
\end{equation}
where $\mu(-u)\in \mathbb R_{++}$ is any positive root of the nonincreasing function:
$$
\vartheta (\mu)=\|{\cal P}_{\mu}(-a)\|_1 -\lambda+s.
$$
\end{lemma}
Like Proposition \ref{prop:mproj-f}, we obtain the following conclusions without proof.
\begin{proposition}\label{Mpprojl1}
Let $z=(A,s)\in {\cal H}$ be a given point with $A=P{\rm Diag}(\lambda (A))P^T$, then
  the following results hold:
\begin{equation}\label{maeq:kfap}
\Pi_{{\cal K}_1}(A,s)=\left
\{
\begin{array}{ll}
(A,s) & \mbox{if } (A,s) \in (A,s) \in {\cal K}_1,\\[2pt]
(0,0) & \mbox{if } (A,s) \in {\cal K}_1^{\circ},\\[2pt]
(P{\rm Diag}\left({\cal P}_{\mu (z)}(\lambda(A))\right)P^T, \mu(z)+s) &  \mbox{otherwise},
\end{array}
\right.
\end{equation}
where $\mu(z)\in \mathbb R_{++}$ is any positive root of the nonincreasing function:
$$
\theta (\mu)=\|{\cal P}_{\mu}(\lambda (A))\|_1 -\mu-s.
$$
\end{proposition}

\begin{proposition}\label{MPprojlinf}
Let $z=(A,s)\in {\cal H}$ be a given point with $A=P{\rm Diag}(\lambda (A))P^T$, then
  the following results hold:
\begin{equation}\label{mmeq:kfap}
\Pi_{{\cal K}_{\infty}}(A,s)=\left
\{
\begin{array}{ll}
(A,s) & \mbox{if } (A,s) \in {\cal K}_{\infty},\\[2pt]
(0,0) & \mbox{if } (a,s) \in {\cal K}_{\infty}^{\circ},\\[2pt]
(A,s)-(P{\rm Diag}\left({\cal P}_{\mu (-z)}(-\lambda (A))\right)P^T, \lambda(-z)-s) &  \mbox{otherwise},
\end{array}
\right.
\end{equation}
where $\mu(-z)\in \mathbb R_{++}$ is any positive root of the nonincreasing function:
$$
\vartheta (\mu)=\|{\cal P}_{\mu}(-\lambda (A))\|_1 -\lambda+s.
$$
\end{proposition}
\subsection{The  Formula for $\Pi_{\K_p}(z)$ when $p \in (1,\infty)$}\label{Sec3.2}
For a given  $z=(A,s)\in {\cal H}$ and any $(u,\mu) \in \mathbb{R}^{m}\times \mathbb{R}$ with $\varphi_p(u) \ne 0$, define $F_p:\mathbb{R}^{m+1}\rightarrow
\mathbb{R}^{m+1}$ by
\begin{equation}\label{eq:FFppp}
F_p(u,\mu;z)=\left[
\begin{array}{l}
\mu \left(\textbf{1}_m^T|u|^p \right)^{-1/p_*}|u|^{p-1}\circ {\rm sgn}(u)+u-\lambda (A)\\[3pt]
\varphi_p(u)-\mu -s
\end{array}
\right
].
\end{equation}
It follows from Proposition \ref{prop:mproj-f}, we have the the following result for $\Pi_{\K_p}(z)$.
\begin{proposition}\label{mmmprop:proj-fap}
Let $z=(A,s)\in {\cal H}$ be  a given point with $A=P{\rm Diag}(\lambda (A))P^T$, then
\begin{equation}\label{mmmeq:kfap}
\Pi_{\K_p}(A,s)=\left
\{
\begin{array}{ll}
(A,s) & \mbox{if } (A,s) \in \K_p,\\[2pt]
(0,0) & \mbox{if } (A,s) \in [\K_p]^{\circ},\\[2pt]
(P{\rm Diag}(u(z))P^T, \mu(z)+s) &  \mbox{otherwise},
\end{array}
\right.
\end{equation}
where $(u(z),\mu(z))\in \mathbb R^m \times \mathbb R_{++}$ is a solution of the system of equations:
$$
F_p(u,\mu;z)=0,
$$
where $F_p$ is defined by (\ref{eq:FFppp}).
\end{proposition}

\subsection{Variational Geometry of $\K_p$ when $p\in [2,\infty)$}\label{Sec3.3}

The second-order derivative of $(\varphi_p\circ\lambda)$ depends on the following proposition about  the spectral function of
a twice differentiable symmetric gauge function $\psi$.

If $p\in [2,\infty)$ and $y \in \mathbb R^m$ with  $y \ne 0$, one has that $\varphi_p$ is twice continuously at $y$. In this case, the gradient of $\varphi_p$ at $y$ is
\begin{equation}\label{eq:grad-phi}
\begin{array}{ll}
\nabla \varphi_p (y)& = \displaystyle \frac{1}{p}\left(\displaystyle \sum_{j=1}^m |y_j|^p \right)^{1/p-1}\left
[
\begin{array}{c}
p |y_1|^{p-1} {\rm sgn}\,(y_1)\\[4pt]
\vdots\\[4pt]
p |y_m|^{p-1} {\rm sgn}\,(y_m)
\end{array}
\right]\\[16pt]
&=\left(\textbf{1}_m^T|y|^p \right)^{-1+1/p}|y|^{p-1}\circ {\rm sgn}(y)\\[12pt]
&=\left(\textbf{1}_m^T|y|^p \right)^{-1/p_*}|y|^{p-1}\circ {\rm sgn}(y).
\end{array}
\end{equation}
The Hessian of  $\varphi_p$ at $y$ is
\begin{equation}\label{eq:H-phi}
\begin{array}{ll}
\nabla^2 \varphi_p (y)& = \left(\displaystyle \frac{1}{p}-1\right)\left(\displaystyle \sum_{j=1}^m |y_j|^p\right)^{1/p-2}\cdot p\left[|y|^{p-1}\circ {\rm sgn}(y)\right]\left[|y|^{p-1}\circ {\rm sgn}(y)\right]^T\\[12pt]
& \quad +\left(\displaystyle \sum_{j=1}^m |y_j|^p\right)^{1/p-1}\cdot (p-1) {\rm diag}\,(|y_1|^{p-2},\ldots, |y_m|^{p-1})\\[16pt]
&= (p-1) |Y|^{p-2}\left(\textbf{1}_m^T|y|^p  \right)^{1/p-1}\\[12pt]
& \quad -(p-1)\left(\textbf{1}_m^T|y|^p  \right)^{1/p-2}\left[|y|^{p-1}\circ {\rm sgn}(y)\right]\left[|y|^{p-1}\circ {\rm sgn}(y)\right]^T\\[12pt]
&= (p-1)\left(\textbf{1}_m^T|y|^p  \right)^{-1/p_*}\left[|Y|^{p-2}-\displaystyle \displaystyle \frac{ |Y|^{p-2}yy^T|Y|^{p-2}}{y^T|Y|^{p-2}y} \right].
\end{array}
\end{equation}

\subsubsection{Tangent sets}
For $(A,s) \in \K_p$, let $A$ have the spectral decomposition $A=P{\rm Diag}w P^T$.  Then we obtain the following formula of the tangent cone:
$$
{\cal T}_{\K_p}(A,s)=\left
\{
\begin{array}{ll}
{\cal H} & \|A\|_p<s,\\[2pt]
\K_p & (A,s)=(0,0),\\[2pt]
\left\{(d_w,d_s):\langle {\rm Diag}\left( |\lambda (A)|^{p-1}\circ {\rm sgn}(\lambda (A))\right),
P^Td_wP \rangle\leq  s^{p/p_*}d_s \right\} &\|A\|_p=s>0.
\end{array}
\right.
$$
For $(A,s) \in [\K_p]^{\circ}$, namely $-s \geq  (\varphi_{p_*}\circ \lambda)(A)$, we have
$$
{\cal T}_{[\K_p]^{\circ}}(A,s)=\left
\{
\begin{array}{ll}
{\cal H} &  \|A\|_{p_*} < -s,\\[2pt]
[\K_p]^{\circ} & (A,s)=(0,0),\\[2pt]
\left\{(d_w,d_s): \langle |\lambda (A)|^{p_*-1}\circ {\rm sgn}(\lambda (A)),
P^Td_w P\rangle\leq -s^{p_*/p}d_s \right\} &\|A\|_{p_*}=-s>0.
\end{array}
\right.
$$
For $(A,s) \in \K$,
$$
{\cal N}_{\K_p}(A,s)=\left
\{
\begin{array}{ll}
\{0\}  & \|A\|_p <s,\\[2pt]
[\K_p]^{\circ} & (A,s)=(0,0),\\[2pt]
\big\{\alpha(P{\rm Diag}(v)P^T, -1):\varphi_{p_*}(v)=1, \langle v, \lambda(A)/s\rangle=1,\alpha\geq 0 \big\} & \|A\|_p=s>0.
\end{array}
\right.
$$
Let $w_{k_1},\ldots, w_{k_r}$ be $r$ distinct values of $m$ eigenvalues of $A$, and $A$ has the spectral decomposition $A=P{\rm Diag}w P^T$ with $w=\lambda(A)$,
namely
$$
w_1=\cdots=w_{k_1} > w_{k_1+1}=\cdots=w_{k_2} > w_{k_2+1}\cdots w_r,
$$
where $k_0=0, k_r=m$. Denote
$$
\alpha_1=\{1,\ldots, k_1\}, \alpha_2=\{k_1+1,\ldots, k_2\}, \ldots, \alpha_r=\{k_{r-1}+1,\ldots, k_r\}.
$$
Define
\begin{equation}\label{eq:nts}
\begin{array}{l}
c_i=\varphi_p(w)^{-p/p_*}|w_i|^{p-1}{\rm sgn}(w_i), b_i=(p-1)\varphi_p(w)^{-p/p_*}|w_i|^{p-2},i=1,\ldots, r, \\[6pt]
a_{ij}=-(p-1)\varphi_p(w)^{1-2p}|w_iw_j|^{p-2}w_iw_j, i,j=1,\ldots, r,\\[6pt]
\delta_{ii}=(p-1)\varphi_p(w)^{-p/p_*}|w_i|^{p-2}+(p-1)\varphi_p(w)^{1-2p}|w_i|^{2(p-1)},i=1,\ldots, r,\\[6pt]
\delta_{ij}=\displaystyle \frac{c_i-c_j}{w_i-w_j},i,j=1,\ldots, r, i\ne j.
\end{array}
\end{equation}
\begin{proposition}\label{lem:m2ndTap}
Let $z \in \K_p$ and $d \in {\cal T}_{ \K_p}(z)$ where $z=(A,s)$ with $A=P{\rm Diag}(\lambda (A))P$, and $d=(d_w,d_s)\in {\cal T}_{{\cal K}_p}(z)$. Then
\begin{equation}\label{eq:2ndTKap}
\begin{array}{l}
{\cal T}^2_{ \K_p}(z,d)\\[20pt]
=\left\{
\begin{array}{lr}
{\cal H} \quad \quad  \quad \quad\quad \quad\quad \quad \quad \quad \quad \quad  \quad \quad \quad \quad \quad\quad \quad\quad \quad \quad \quad   d \in \mbox{int }{\cal T}_{ \K_p}(z),\\[6pt]
{\cal T}_{ \K_p}(z) \quad \quad  \quad \quad\quad \quad\quad \quad \quad \quad\quad \quad  \quad \quad\quad \quad\quad \quad \quad \quad  \quad \quad \quad \quad \quad \quad z=0,\\[6pt]
\left\{(\xi_w,\xi_s):
\begin{array}{l}
s^{-p/p_*}\langle {\rm Diag}\left(|\lambda (A)|^{p-1}\circ {\rm sgn}(\lambda (A))\right),
P^T\xi_wP \rangle \\[8pt]
 +\displaystyle \sum_{i,j=1}^r a_{ij}{\rm Tr}
(P_{\alpha_i}^Td_AP_{\alpha_i}){\rm Tr}
(P_{\alpha_j}^Td_AP_{\alpha_j})\\[8pt]
\quad \, +\displaystyle \sum_{i,j=1}^r \delta_{ij}\displaystyle \sum_{i'\in \alpha_i,j'\in \alpha_j}(p_{i'}^Td_Ap_{j'})^2-\xi_s\leq 0
\end{array}
\right\} \,\, \mbox{otherwise},
\end{array}
\right.
\end{array}
\end{equation}
where $a_{ij}$ and $\delta_{ij}$ are defined by (\ref{eq:nts}).
\end{proposition}

\subsubsection{Directional derivative and B-subdifferential of $\Pi_{\K_p}(z)$}
Let $(u(z),\mu(z))\in \mathbb R^m \times \mathbb R_{++}$ be a solution of the system of equations:
$$
F_p(u,\mu;z)=0,
$$
where $F_p$ is defined by (\ref{eq:FFppp}). Let $w_1(z):=u_{k_1(z)}(z),\ldots, w_{r(z)}(z):=u_{k_r(z)}(z)$ be $r(z)$ distinct values of $m$ values  of $u(z)$, and $X(z)$ has the spectral decomposition $X(z)=P{\rm Diag}u(z) P^T$ with $u(z)=\lambda(X(z))$,
namely
$$
u_1(z)=\cdots=u_{k_1(z)}(z) > u_{k_1(z)+1}(z)=\cdots=u_{k_2(z)}(z) > u_{k_2(z)+1}(z)\cdots u_{k_r(z)}(z),
$$
where $k_0=0, k_r(z)=m$. Denote
$$
\alpha_1(z)=\{1,\ldots, k_1(z)\}, \alpha_2(z)=\{k_1(z)+1,\ldots, k_2(z)\}, \ldots, \alpha_r(z)=\{k_{r-1}(z)+1,\ldots, k_r(z)\}.
$$
Define
\begin{equation}\label{eq:ntszz}
\begin{array}{l}
c_i(z)=\varphi_p(u(z))^{-p/p_*}|u_i(z)|^{p-1}{\rm sgn}(u_i(z)), i=1,\ldots, r(z), \\[6pt]
b_i(z)=(p-1)\varphi_p(u(z))^{-p/p_*}|u_i(z)|^{p-2},i=1,\ldots, r(z), \\[6pt]
a_{ij(z)}=-(p-1)\varphi_p(u(z))^{1-2p}|u_i(z)u_j(z)|^{p-2}u_i(z)u_j(z), i,j=1,\ldots, r(z),\\[6pt]
\delta_{ii}(z)=(p-1)\varphi_p(u(z))^{-p/p_*}|u_i(z)|^{p-2}\\[6pt]
 \quad \quad \quad \quad +(p-1)\varphi_p(u(z))^{1-2p}|u_i(z)|^{2(p-1)},i=1,\ldots, r(z),\\[6pt]
\delta_{ij}(z)=\displaystyle \frac{c_i(z)-c_j(z)}{w_i(z)-w_j(z)},i,j=1,\ldots, r(z), i\ne j,\\[12pt]
\pi_{ij}(z)=(1+\mu(z)\delta_{ij}(z))^{-1}, i,j=1,\ldots, r(z).
\end{array}
\end{equation}
\begin{theorem}\label{HelppthmthPiDirction}
 Let $z=(A,s)\in {\cal H}$ be a given point with $A=P{\rm Diag}(\lambda (A))P^T$, then the following results hold:
\begin{itemize}
\item[{\rm (i)}] If $z=(A,s) \in {\rm int } \K_p$, then $\Pi'_{\K_p}(z;d_z)=d_z$;
\item[{\rm (ii)}] If $z=(A,s) \in {\rm int } [\K_p]^{\circ}$, then $\Pi'_{\K_p}(z;d_z)=0$;
\item[{\rm (iii)}] If $z=(A,s) \in {\rm int }[\mathbb S^{m}\times \mathbb R\setminus (\K_p\cup[\K_p]^{\circ})]$, then
\begin{equation}\label{meq:projD-1pg2}
\begin{array}{l}
\Pi'_{\K_p}(z;d_z)=
\left[
\begin{array}{c}
P\left[ {\rm Diag}\left(M_p(z)^{-1}{\rm diag}(P^Td_AP)\right)+{\cal B}(z)\circ  P^Td_AP\right]P^T
\\[14pt]
d_s
\end{array}
\right] \\[4pt]
\quad -\displaystyle \frac{[\langle P {\rm Diag}\left[M_p(z)^{-1}\nabla \varphi_p(u(z))\right]P^T, d_A \rangle-d_s]}{1+\nabla \varphi_p(u(z))^TM_p(z)^{-1}\varphi_p(u(z))}\left[
\begin{array}{c}
P {\rm Diag}\left[M_p(z)^{-1}\nabla \varphi_p(u(z))\right]P^T\\[14pt]
-1
\end{array}
\right],
\end{array}
\end{equation}
where ${\cal B}(z)$ is defined by (\ref{2ndpsimuB}),
 $X(z)=P{\rm Diag}(u(z))P^T$, $M_p(z)=I_m+\mu(z)\nabla^2\varphi_p(u(z))$ with $(u(z),\mu(z))\in \mathbb R^m \times \mathbb R_{++}$ being a solution of the system of equations:
$$
F_p(u,\mu;z)=0.
$$
\item[{\rm (iv)}] If $z=(A,s) \in {\rm bdry} \K_p \setminus\{0\}$, then
  \begin{equation}\label{meq:projD-2}
\begin{array}{l}
\Pi'_{\K_p}(z;d_z)=
\left[
\begin{array}{c}
d_A\\[14pt]
d_s
\end{array}
\right]-\displaystyle \frac{[\nabla (\varphi_p\circ \lambda)(A)^*d_A -d_s]_+}{1+\|\nabla (\varphi_p\circ \lambda)(A)\|^2_F}\left[
\begin{array}{c}
\nabla (\varphi_p\circ \lambda)(A)\\[14pt]
-1
\end{array}
\right],
\end{array}
\end{equation}
where
$$
\nabla (\varphi_p\circ \lambda)(A)=\|A\|_p^{-p/p_*}\langle P{\rm Diag}(|\lambda (A)|^{p-1}\circ {\rm sgn}(\lambda (A)))P^T.
$$
\item[{\rm (v)}] If $z=(A,s) \in {\rm bdry} [\K]^{\circ} \setminus\{0\}$, then
  \begin{equation}\label{eq:projD-3}
\begin{array}{l}
\Pi'_{\K}(z;d_z)=
\displaystyle \frac{[\nabla (\varphi_{p_*}\circ \lambda)(A)^*d_A +d_s]_+}{1+\|\nabla (\varphi_{p_*}\circ \lambda)(A)\|^2_F}\left[
\begin{array}{c}
\nabla ((\varphi_{p_*}\circ \lambda)(A)\\[14pt]
1
\end{array}
\right],
\end{array}
\end{equation}
where
$$
\nabla (\varphi_{p_*}\circ \lambda)(A)=\|A\|_{p_*}^{-p_*/p}\langle P{\rm Diag}(|\lambda (A)|^{p_*-1}\circ {\rm sgn}(\lambda (A)))P^T.
$$
\item[{\rm (vi)}] If $z=(A,s)=(0,0)$, then $\Pi'_{\K_p}(z;d_z)=\Pi_{\K_p}({\rm d}_z)$.
\end{itemize}
\end{theorem}

\subsection{Variational Geometry of $\K_p$ when $p\in (1,2)$}\label{Sec3.4}
If $p\in (1,2)$, then $p_* \in (2, \infty)$. Then for $y \in \mathbb R^m$ with  $y \ne 0$, one has that $\varphi_{p_*}$ is twice continuously at $y$. In this case, the gradient of $\varphi_{p_*}$ at $y$ is
\begin{equation}\label{eq:grad-phida}
\nabla \varphi_{p_*} (y)=\left(\textbf{1}_m^T|y|^{p_*} \right)^{-1/p}|y|^{p_*-1}\circ {\rm sgn}(y)
\end{equation}
and the Hessian of $\varphi_{p_*}$ at $y$ is
\begin{equation}\label{eq:H-phida}
\nabla^2 \varphi_{p_*} (y)=(p_*-1)\left(\textbf{1}_m^T|y|^{p_*}  \right)^{1/p_*-1}\left[|Y|^{p_*-2}-\displaystyle \displaystyle \frac{ |Y|^{p_*-2}yy^T|Y|^{p_*-2}}{y^T|Y|^{p_*-2}y} \right].
\end{equation}
\subsubsection{Tangent Sets}
When $p \in (1,2)$, we also have for $(A,s) \in \K_p$ that
$$
{\cal T}_{\K_p}(A,s)=\left
\{
\begin{array}{ll}
\mathbb{S}^{m}\times\mathbb R &  \|A\|_p<s,\\[2pt]
\K_p & (A,s)=(0,0),\\[2pt]
\left\{(d_A,d_s): \|A\|_p^{-p/p_*}\langle P{\rm Diag}(|\lambda (A)|^{p-1}\circ {\rm sgn}(\lambda (A)))P^T,
d_A \rangle\leq d_s \right\} & \|A\|_p=s>0.
\end{array}
\right.
$$
For $(A,s) \in [\K_p]^{\circ}$, namely $-s \geq \|A\|_{p_*}$, we have
$$
{\cal T}_{[\K_p]^{\circ}}(A,s)=\left
\{
\begin{array}{ll}
\mathbb{S}^{m}\times\mathbb R &  \|A\|_{p_*} < -s,\\[2pt]
[\K_p]^{\circ} & (A,s)=(0,0),\\[2pt]
\left\{(d_A,d_s): \|A\|_{p_*}^{-p_*/p}\langle P{\rm Diag}(|\lambda (A)|^{p_*-1}\circ {\rm sgn}(\lambda (A)))P^T,
d_A \rangle\leq -d_s \right\} & \|A\|_{p_*}=-s>0.
\end{array}
\right.
$$

\subsubsection{Directional derivative and B-subdifferential of $\Pi_{\K_p}(z)$}
For a given  $z=(A,s)\in \mathbb{R}^{m}\times \mathbb{R}$ and any $(u,\mu) \in \mathbb{R}^{m}\times \mathbb{R}$ with $\|u\|_p \ne 0$, define $F_{p_*}:\mathbb{R}^{m+1}\rightarrow
\mathbb{R}^{m+1}$ by
\begin{equation}\label{eq:Fpstar}
F_{p_*}(u,\mu;z)=\left[
\begin{array}{l}
\mu \left(\textbf{1}_m^T|u|^{p_*} \right)^{-1/p}|u|^{p_*-1}\circ {\rm sgn}(u)+u-\lambda(A)\\[3pt]
\|u\|_{p_*}-\mu -s
\end{array}
\right
].
\end{equation}
Define
\begin{equation}\label{eq:ntszzstar}
\begin{array}{l}
\widehat c_i(z)= \varphi_{p_*}(u(z))^{-p_*/p}|u_i(z)|^{p^*-1}{\rm sgn}(u_i(z)), i=1,\ldots, r(z), \\[6pt]
\widehat b_i(z)= (p_*-1)\varphi_{p_*}(u(z))^{-p_*/p}|u_i(z)|^{p_*-2},i=1,\ldots, r(z), \\[6pt]
\widehat a_{ij(z)}=-(p_*-1)\varphi_{p_*}(u(z))^{1-2p_*}|u_i(z)u_j(z)|^{p_*-2}u_i(z)u_j(z), i,j=1,\ldots, r(z),\\[6pt]
\delta_{ii}(z)=(p_*-1)\varphi_{p_*}(u(z))^{-p_*/p}|u_i(z)|^{p_*-2}\\[6pt]
\quad \quad  \quad +(p_*-1)\varphi_{p_*}(u(z))^{1-2p_*}|u_i(z)|^{2(p_*-1)},i=1,\ldots, r(z),\\[6pt]
\widehat \delta_{ij}(z)=\displaystyle \frac{\widehat c_i(z)-\widehat c_j(z)}{w_i(z)-w_j(z)},i,j=1,\ldots, r(z), i\ne j,\\[12pt]
\widehat \pi_{ij}(z)=(1+\mu(z)\widehat\delta_{ij}(z))^{-1}, i,j=1,\ldots, r(z),\\[12pt]
M_{p_*}(z)=I_m+\mu(z)\nabla^2\varphi_{p_*}(u(z))
\end{array}
\end{equation}
and
\begin{equation}\label{2ndpsimuBstar}
\begin{array}{ll}
{\cal B}_*(z)=& -\left[
\begin{array}{ccc}
\widehat \pi_{11}(z)I_{|\alpha_1(z)|}  & &\\[4pt]
&  \ddots  &\\[4pt]
&  & \widehat \pi_{r(z)r(z)}(z)I_{|\alpha_r(z)|}
\end{array}
\right]\\[18pt]
& +\left [
\begin{array}{ccc}
\widehat \pi_{11}(z)\textbf{1}_{|\alpha_1(z)|}\textbf{1}_{|\alpha_1(z)|}^T &\cdots & \widehat \pi_{1r(z)}(z)\textbf{1}_{|\alpha_1(z)|}\textbf{1}_{|\alpha_{r(z)}(z)|}^T\\[10pt]
 \vdots &\vdots& \vdots\\[10pt]
\widehat \pi_{r(z)1}(z)\textbf{1}_{|\alpha_{r(z)}(z)|}\textbf{1}_{|\alpha_1(z)|}^T & \cdots &\widehat \pi_{r(z)r(z)}(z)\textbf{1}_{|\alpha_{r(z)}(z)|}\textbf{1}_{|\alpha_{r(z)}(z)|}^T
\end{array}
\right].
\end{array}
\end{equation}
Using that $z=\Pi_{\K_p}(z)+\Pi_{[\K_p]^{\circ}}(z)$, we obtain the following result.
\begin{theorem}\label{HelppthmthPiDirction}
 Let $z=(A,s)\in {\cal H}$ be a given point with $A=P{\rm Diag}(\lambda (A))P^T$, then the following results hold:
\begin{itemize}
\item[{\rm (i)}] If $z=(A,s) \in {\rm int } [\K_p]^{\circ}$, then $\Pi'_{\K_p}(z;d_z)=0$;
\item[{\rm (ii)}] If $z=(A,s) \in {\rm int }\K_p$, then $\Pi'_{\K_p}(z;d_z)=d_z$;
\item[{\rm (iii)}] If $z=(A,s) \in {\rm int }[\mathbb S^{m}\times \mathbb R\setminus (\K_p\cup[\K_p]^{\circ})]$, then
\begin{equation}\label{meq:projD-1pg2}
\begin{array}{l}
\Pi'_{\K_p}(z;d_z)=
\left[
\begin{array}{c}
d_A-P\left[ {\rm Diag}\left(M_{p_*}(z)^{-1}{\rm diag}(P^Td_AP)\right)+{\cal B}_*(z)\circ  P^Td_AP\right]P^T\\[14pt]
0
\end{array}
\right] \\[4pt]
\quad +\displaystyle \frac{[\langle P {\rm Diag}\left[M_{p_*}(z)^{-1}\nabla \varphi_{p_*}(u(z))\right]P^T, d_A\rangle +d_s]}{1+\nabla \varphi_{p_*} (u(z))^TM_{p_*}(z)^{-1}\nabla \varphi_{p_*} (u(z))}\left[
\begin{array}{c}
P {\rm Diag}\left[M_{p_*}(z)^{-1}\nabla \varphi_{p_*}(u(z))\right]P^T\\[14pt]
1
\end{array}
\right],
\end{array}
\end{equation}
where $X(z)=P{\rm Diag}(u(z))P^T$ and $(u(z),\mu(z))\in \mathbb R^m \times \mathbb R_{++}$ is a solution of the system of equations:
$$
F_{p_*}(u,\mu;z)=0.
$$
\item[{\rm (iv)}] If $z=(A,s) \in {\rm bdry} \K_p \setminus\{0\}$, then
  \begin{equation}\label{meq:projD-2}
\begin{array}{l}
\Pi'_{\K_p}(z;d_z)=
\left[
\begin{array}{c}
d_A\\[14pt]
d_s
\end{array}
\right]-\displaystyle \frac{[\nabla (\varphi_p\circ \lambda)(A)^*d_A -d_s]_+}{1+\|\nabla (\varphi_p\circ \lambda)(A)\|^2_F}\left[
\begin{array}{c}
\nabla (\varphi_p\circ \lambda)(A)\\[14pt]
-1
\end{array}
\right],
\end{array}
\end{equation}
where
$$
\nabla (\varphi_p\circ \lambda)(A)=\|A\|_p^{-p/p_*}\langle P{\rm Diag}(|\lambda (A)|^{p-1}\circ {\rm sgn}(\lambda (A)))P^T.
$$
\item[{\rm (v)}] If $z=(A,s) \in {\rm bdry} [\K]^{\circ} \setminus\{0\}$, then
  \begin{equation}\label{eq:projD-3}
\begin{array}{l}
\Pi'_{\K}(z;d_z)=
\displaystyle \frac{[\nabla (\varphi_{p_*}\circ \lambda)(A)^*d_A +d_s]_+}{1+\|\nabla (\varphi_{p_*}\circ \lambda)(A)\|^2_F}\left[
\begin{array}{c}
\nabla ((\varphi_{p_*}\circ \lambda)(A)\\[14pt]
1
\end{array}
\right],
\end{array}
\end{equation}
where
$$
\nabla (\varphi_{p_*}\circ \lambda)(A)=\|A\|_{p_*}^{-p_*/p}\langle P{\rm Diag}(|\lambda (A)|^{p_*-1}\circ {\rm sgn}(\lambda (A)))P^T.
$$
\item[{\rm (vi)}] If $z=(A,s)=(0,0)$, then $\Pi'_{\K_p}(z;d_z)=\Pi_{\K_p}({\rm d}_z)$.
\end{itemize}
\end{theorem}

\section{{Conclusion}} \label{Section-Conclusion}
This paper developed variational analysis of orthogonally invariant  norm cone of symmetric matrices including  the computation of the tangent cone, the  normal cone, and the
second-order tangent set of an orthogonally invariant  norm cone, as well as the formulas for directional derivative and B-subdifferential of the projection operator onto a norm cone.  These set of general results are applied to Schatten $p$-norm cone for $p \in (1,\infty)$, particularly to the second-order cone.  There are many  topics  left to study, for instances
 the optimality theories of Problem (\ref{Pconic})  are needed to be established, and
   the perturbation analysis of Problem (\ref{Pconic}) is an important topic worth studying.
   Noting that the orthogonally invariant  norm is a spectral function of symmetric matrices,
    Cui,  Ding and Zhao (2017) \cite{CZhao2017} already discussed quadratic growth conditions for convex matrix optimization problems associated with spectral functions. Therefore, there are many theoretical problems to study in orthogonally invariant  norm conic
    optimization.


\end{document}